\documentclass[10pt]{article}
\usepackage{amscd,latexsym,amsmath,amssymb,amsfonts,graphicx,fancyhdr,setspace,xypic}
\usepackage{latexsym,amsmath,amssymb,amsfonts,verbatim,hyperref,chngcntr}
\usepackage{amsmath,amscd}

\author{
Junwu Tu\thanks{Mathematics Department, University of Oregon, Eugene OR 97403, USA, {\em e-mail:}{\tt junwut@uoregon.edu}}}
\title{ Homological mirror symmetry and Fourier-Mukai transform }
\date{}

\DeclareFontFamily{U}{rsf}{}
\DeclareFontShape{U}{rsf}{m}{n}{
  <5> <6> rsfs5 <7> <8> <9> rsfs7 <10-> rsfs10}{}
\DeclareMathAlphabet{\mathscr}{U}{rsf}{m}{n}

\DeclareMathAlphabet{\mathgth}{U}{euf}{m}{n}

\DeclareFontFamily{U}{cyr}{}
\DeclareFontShape{U}{cyr}{m}{n}{
  <5> wncyr5 <6> wncyr6 <7> wncyr7 <8> wncyr8 <9> wncyr9 <10-> wncyr10}{}
\DeclareMathAlphabet{\mathcyr}{U}{cyr}{m}{n}

\makeatletter
\def\operator@font{\sf}
\makeatother

\newcommand{\cA}{{\mathscr A}}

\newcommand{\cW}{{\mathscr W}}

\newcommand{\sL}{{\mathcal L}}

\newcommand{\ccH}{{\mathscr H}}

\newcommand{\cM}{{\mathscr M}}
\newcommand{\cN}{{\mathscr N}}
\newcommand{\cO}{{\mathscr O}}
\newcommand{\cP}{{\mathscr P}}

\newcommand{\even}{{\mathsf{even}}}

\newcommand{\tw}{{\mathsf{tw}}}

\newcommand{\can}{{\mathsf{can}}}
\newcommand{\dR}{{\mathsf{dR}}}
\newcommand{\Fuk}{\mathsf{Fuk}}
\newcommand{\MC}{\mathsf{MC}}

\newcommand{\D}{{\mathsf D}_{\mathsf{coh}}^b}
\newcommand{\sing}{{\mathsf{sing}}}

\newcommand{\chk}{{\scriptscriptstyle\vee}}
\newcommand{\R}{\mathbb{R}}

\newcommand{\db}{{\overline{\partial}}}
\newcommand{\hol}{{\mathsf{hol}}}

\newcommand{\val}{{\mathsf{val}}}

\newcommand{\interior}{{\mathsf{int}}}

\newcommand{\bone}{{\mathbf 1}}

\DeclareMathOperator{\End}{End}
\DeclareMathOperator{\Hom}{Hom}

\DeclareMathOperator{\Tw}{Tw}

\DeclareMathOperator{\id}{id}
\DeclareMathOperator{\ev}{ev}

\DeclareMathOperator{\Ext}{Ext}

\DeclareMathOperator{\sym}{sym}

\newcommand{\ra}{\rightarrow}

\newcommand{\C}{\mathbb{C}}

\newcommand{\Z}{\mathbb{Z}}

\newcommand{\remark}{\noindent\textbf{Remark: }}
\newcommand{\proof}{\noindent\textbf{Proof. }}

\newtheorem{theorem}{Theorem}[section]
\newtheorem{assumption}[theorem]{Assumption}

 \newtheorem{lemma}[theorem]{Lemma}
 
 \newtheorem{proposition}[theorem]{Proposition}

 \newtheorem{definition}[theorem]{Definition}
 \newtheorem{definition-theorem}[theorem]{Definition-Theorem}

\renewcommand{\phi}{\varphi}

\makeatletter
\newcommand*{\doublerightarrow}[2]{\mathrel{
  \settowidth{\@tempdima}{$\scriptstyle#1$}
  \settowidth{\@tempdimb}{$\scriptstyle#2$}
  \ifdim\@tempdimb>\@tempdima \@tempdima=\@tempdimb\fi
  \mathop{\vcenter{
    \offinterlineskip\ialign{\hbox to\dimexpr\@tempdima+1em{##}\cr
    \rightarrowfill\cr\noalign{\kern.01ex}
    \rightarrowfill\cr}}}\limits^{\!#1}_{\!#2}}}
\newcommand*{\triplerightarrow}[1]{\mathrel{
  \settowidth{\@tempdima}{$\scriptstyle#1$}
  \mathop{\vcenter{
    \offinterlineskip\ialign{\hbox to\dimexpr\@tempdima+1em{##}\cr
    \rightarrowfill\cr\noalign{\kern.01ex}
    \rightarrowfill\cr\noalign{\kern.01ex}
    \rightarrowfill\cr}}}\limits^{\!#1}}}

\newcommand{\Rnum}[1]{\expandafter\@slowromancap\romannumeral #1@}
\makeatother
\numberwithin{equation}{section}
\counterwithout{paragraph}{subsection}
\counterwithin{paragraph}{section}
\begin{document}
\maketitle
\begin{abstract}

We interpret symplectic geometry as certain sheaf theory by constructing a sheaf of curved $A_\infty$ algebras which in some sense plays the role of a ``structure sheaf" for symplectic manifolds. An interesting feature of this ``structure sheaf" is that the symplectic form itself is part of its curvature term. Using this interpretation homological mirror symmetry can be understood by well-known duality theories in mathematics: Koszul duality or Fourier-Mukai transform. In this paper we perform the above constructions over a small open subset inside the smooth locus of a Lagrangian torus fibration. In a subsequent work we shall use the language of derived geometry to obtain a global theory over the whole smooth locus. However we do not know how to extend this construction to the singular locus. As an application of the local theory we prove a version of homological mirror symmetry between a toric symplectic manifold and its Landau-Ginzburg mirror. 


\end{abstract}
\section{Introduction}

\paragraph{Backgrounds and histories.} Homological mirror symmetry conjecture was proposed by M. Kontsevich~\cite{Kont} in an address to the 1994 International Congress of Mathematicians, aiming to give a mathematical framework to understand the mirror phenomenon originated from physics. Roughly speaking this conjecture predicts a quasi-equivalence between two $A_\infty$ triangulated categories $\Fuk(M)$ and $\D(M^\chk)$ naturally associated to a symplectic manifold $M$ and a complex manifold $M^\chk$. Note that despite of the notation, the mirror manifold $M^\chk$ is not uniquely determined by $M$. After nearly two decades since Kontsevich's proposal, his conjecture  has been generalized and has been proven in a lot of deep and inspiring situations~\cite{PZ},~\cite{KS},~\cite{Seidel},~\cite{Sh},~\cite{Seidel2},~\cite{AKO},~\cite{FLTZ}, and definitely many more~\footnote{It seems to the author impossible to make a complete list, so he apologize for this.}.

In spite of our increasing knowledge of this conjecture, less is understood about the mathematical reason behind it. The first attempt towards a general mathematical theory to understand the mirror phenomenon is given by A. Strominger, S-T. Yau and E. Zaslow. In~\cite{SYZ} they proposed a geometric picture to produce mirror pairs: mirror duality should arise between (special) Lagrangian torus fibrations and the associated dual fibrations. Indeed the SYZ proposal was successfully realized in the semiflat cases, showing beautiful interactions between special Lagrangians and stable holomorphic vector bundles~\cite{LYZ},~\cite{AP},~\cite{CL}. In the toric Calabi-Yau case, see~\cite{CLL}. To get interesting symplectic manifolds and complex manifolds (other than torus) one has to allow these fibrations to have singular fibers for which we refer to the massive work of M. Gross and B. Siebert on toric degenerations~\cite{GS},~\cite{GS2}. 

While this gives a nice theory to produce mirror pairs, it gives less hint on why mirror pairs produced from SYZ proposal should interchange symplectic geometry with complex geometry~\footnote{The original SYZ explanation for this was from the so-called T-duality in physics.}. The main purpose of the current paper is to suggest an answer to this question. Namely we show the mirror duality between symplectic geometry of a Lagrangian torus fibration and complex geometry of its dual fibration is in fact a well-known duality in mathematics: Koszul duality or its global version Fourier-Mukai transform.

The current paper grew out of understanding an algebraic framework for K. Fukaya, Y.-G. Oh, H. Ohta and K. Ono's series of papers~\cite{FOOOtoric},\cite{FOOOtoric2},\cite{FOOOtoric3}. There are also inspiring works of Fukaya~\cite{Fukaya2} and Seidel~\cite{Seidel3}.

\paragraph{An example.} We begin with a simple example illustrating the main ideas. Let $\R$ be endowed with a linear coordinate $x$, and let $\R^\chk$ be endowed with the dual coordinate $y^\chk$ (this choice of notation will be clear later). The cotangent bundle $T^*(\R)=\R\times \R^\chk$ has a canonical symplectic form $\omega:=dx\wedge dy^\chk$. To have a Lagrangian torus fibration we consider the quotient space $M:=\R \times (\R^\chk/\Z^\chk)$. Since $\omega$ is translation invariant, it descends to a symplectic form on $M$. The projection map  $\pi: M\ra \R$ onto the first component defines a Lagrangian torus fibration.


Consider a complex vector bundle over the base $\R$ whose fiber over a point $u\in \R$ is the cohomology group $H^*(\pi^{-1}(u),\C)$. Denote by $\ccH$ the corresponding sheaf of $C^\infty$-sections. It is well-known that the sheaf $\ccH$ is in fact a D-module endowed with the Gauss-Manin connection $\nabla$. Moreover the cup product on cohomology defines an algebra structure on $\ccH$ which is $\nabla$-flat. Thus its de Rham complex $\Omega_\R^*(\ccH)$ is a differential graded algebra. We remark that the construction of this differential graded algebra does not involve the symplectic structure on $M$.

To encode the symplectic structure into the algebra $\Omega_\R^*(\ccH)$, we can add a curvature term to it which is just the symplectic form $\omega$ itself (up to sign)! More precisely if we denote by $e=dy^\chk$ the $\nabla$-flat integral generator for $\ccH^1$ (degree one part), then $-\omega$ can be viewed as an element of $\Omega_\R^*(\ccH)$ by writing it as $e\otimes dx$. Note that this element $-\omega$ is of even degree, and it is closed in $\Omega_\R^*(\ccH)$. Since the differential graded algebra $\Omega^*(\ccH)$ is supercommutative any such element can be viewed as a curvature term. Thus we have obtained a sheaf of curved differential graded algebras over the base manifold $\R$ which we denote by $\cO^{\omega,\can}$\footnote{The notion is due to the usage of the sheaf $\ccH$ which is fiber-wise cohomology ring which was referred to as the canonical model in~\cite{FOOO2}.} in the following. We may think of it as a structure sheaf of symplectic geometry of $M$.

\paragraph{From symplectic functions to holomorphic functions: Koszul duality.} To understand the mirror phenomenon which interchanges symplectic geometry with complex geometry, let us compute the Koszul dual algebra of $\cO^{\omega,\can}$. We shall work over the base ring $\Omega^*_{\R}$, the complex valued de Rham complex  of the base manifold $\R$. If we ignore the curvature term $-\omega=e\otimes dx$, the underlying differential graded algebra structure of $\cO^{\omega,\can}$ is simply the exterior algebra generated by $e$ in degree one over $\Omega^*_{\R^\chk}$. Hence its Koszul dual algebra is the symmetric algebra $\sym_{\Omega^*_{\R}}(y)$ generated by a variable $y:=e^\chk[-1]$ which is of degree zero (since $e$ has degree one).

For the curvature term observe that $-\omega=e\otimes dx$ is a linear curvature since we are working over $\Omega^*_{\R}$. It is well-known in Koszul dual theory how to deal with linear curvatures: they introduce additional differential in Koszul dual algebras. The resulting Koszul dual algebra of $\cO^{\omega,\can}$ is thus the symmetric algebra $\sym_{\Omega^*_{\R}}(y)\cong \sym(y)\otimes \Omega^*_{\R}$ endowed with a differential acting on it by the operator $(\partial_{x}+\partial_{y}) dx$. If we perform a change of variable by
\[ x\mapsto x, \;\;  y\mapsto -\sqrt{-1} y, \;\; dx \mapsto d\overline{z};\]
this Koszul dual algebra is identified with the Dolbeault algebra resolving the structure sheaf of the complex manifold $T\R=\R\times \R$ with its canonical complex structure\footnote{This assertion is almost correct except we are only using polynomial function in the $y$ variable as opposed to smooth functions. In the main body of the paper we shall use smooth functions.}. From this example we see that the algebraic reason that Koszul duality interchanges a linear curvature term with an additional differential in its Koszul dual, when applied to this situation, gives a direct link between the symplectic structure and its mirror complex structure.

\paragraph{Approaching HMS via duality of algebras.} From the previous example we see that the sheaf of symplectic functions is ``Koszul dual" to the sheaf of holomorphic functions.  In general we propose to understand the homological mirror symmetry conjecture in the following steps:
\begin{itemize}
\item[\Rnum{1}.] Associated to any Lagrangian torus fibration $\pi: M\ra B$ construct a sheaf of $A_\infty$ algebras $\cO^\omega$ which plays the role of a structure sheaf in symplectic geometry;
\item[\Rnum{2}.] Construct another sheaf $\cO^\hol$ over $B$, which in some sense is ``Koszul dual" to $\cO^\omega$;
\item[\Rnum{3}.] Associated to any Lagrangians (with unitary local systems) in $M$ construct a module over $\cO^\omega$, and show that the Fukaya category $\Fuk(M)$ fully faithfully embeds into the category of modules over $\cO^\omega$;
\item[\Rnum{4}.] Understand the module correspondence between the ``Koszul dual" algebras $\cO^\omega$ and $\cO^\hol$.
\end{itemize}

\noindent 
Here the notion of a sheaf of $A_\infty$ algebras used in $\Rnum{1}$ needs more explanation: locally over small open subsets in $B^\interior$ we can construct honest sheaves of $A_\infty$ algebras, i.e. transition functions for gluing the underlying sheaves are in fact linear $A_\infty$ homomorphisms. Globally we need to take care of possible wall-crossing discontinuity by gluing these local sheaves of $A_\infty$ algebras up to coherent homotopy on intersections, triple intersections, quadruple intersections, and so forth, producing a \emph{homotopy sheaf} of $A_\infty$ algebras.

The current paper contains partial results in all four steps mentioned above. Here ``partial" mainly means that we work out the constructions over a small local open subset inside the smooth locus $B_0$ of a Lagrangian torus fibration. The global constructions of $\cO^\omega$ and its mirror complex manifold over the whole smooth locus $B_0$ will appear in forthcoming works~\cite{Tu},~\cite{Tu2}. We should also confess that it is not clear at present how the singular locus $B^\sing$ should enter into this study.

There are lots of advantages in this new approach to homological mirror symmetry conjecture. The first advantage is a natural construction of  (local) mirror functors (see Sections~\ref{sec:koszul} and~\ref{sec:fourier}) using classical Koszul duality theory of modules. Note that functors constructed in this way are $A_\infty$ functors with explicit formulas. It is also conceptually clearer in this approach how lots of ad-hoc constructions in mirror symmetry should enter into the theory. For instance the inclusion of a B-field to complexify symplectic moduli, or the appearance of quantum corrections in deforming the semi-flat complex structure in mirror complex manifold constructions. Finally we believe this approach offers a way to prove homological mirror symmetry conjecture in an abstract form (i.e. without computing both sides explicitly).

In the remaining part of the introduction we give an overview of materials contained in each section.
\paragraph{Section~\ref{sec:symp}.} In this section we deal with step $\Rnum{1}$ over a small open subset $U$ inside the smooth locus $B_0$ of a Lagrangian torus fibration $\pi: M\ra B$. Denote again by $\pi: M(U)\ra U$ the projection map. Consider the sheaf $\ccH$ of $C^\infty$ sections of $R\pi_*\Lambda^\pi\otimes C^\infty_U$ where $\Lambda^\pi$ is certain relative Novikov ring (see Section~\ref{sec:symp} for its precise definition). The Bott-Morse Lagrangian Floer theory developed in~\cite{FOOO} and~\cite{Fukaya} endows for each Lagrangian fiber $L_u$ an $A_\infty$ algebra structure on $H^*(L_u,\Lambda^\pi)$. We shall call this $A_\infty$ algebra Fukaya algebra of $L_u$. For our purpose we need to consider this $A_\infty$ structure as a family over $U$. This construction has been taken care of by K. Fukaya in~\cite{Fukaya} which we follow in this paper.

Just as in the one dimensional example, the sheaf $\ccH$ is a D-module with the Gauss-Manin connection. But this time the structure maps $m_k$ are not $\nabla$-flat in general.  Our main observation is the following compatibility between the D-module structure on $\ccH$ and its $A_\infty$ structure.
\[[\nabla, m_k]=\sum_{i=1}^{k+1} m_{k+1}(\cdots,\omega,\cdots)\]
This equation is a simple consequence of cyclic symmetry proved in~\cite{Fukaya}. An immediate corollary of these compatibility equations is the following theorem.

\begin{theorem}
There is a curved $A_\infty$ algebra structure on the de Rham complex $\Omega^*_U(\ccH)$ whose curvature is given by $m_0-\omega$.
\end{theorem}

\medskip
\noindent We denote this sheaf of curved $A_\infty$ algebras over $U$ by $\cO^{\omega,\can}_{M(U)}$ (or simply $\cO^{\omega,\can}$). The notation is due to the usage of the $A_\infty$ algebra $H^*(L_u,\Lambda^\pi)$ on each fiber of $\ccH$ which was called canonical models in~\cite{FOOO2}. If we replace the canonical model $H^*(L_u,\Lambda^\pi)$ by the full de Rham complex $\Omega^*(L_u,\Lambda^\pi)$ on each fiber $L_u$, then same construction yields an $A_\infty$ structure on the relative de Rham complex of $\pi$. We shall denote this sheaf of $A_\infty$ algebras by $\cO^\omega_{M(U)}$ (or simply $\cO^\omega$). There is a subtle difference between the two sheaves of $A_\infty$ algebras $\cO^\omega$ and $\cO^{\omega,\can}$ in view of mirror symmetry which is explained in Sections~\ref{sec:koszul} and~\ref{sec:fourier}. Roughly speaking the ``mirror" of the $\cO^{\omega,\can}$ is the tangent bundle $TU$ while that of $\cO^{\omega}$ is the dual torus bundle over $U$.

\paragraph{Section~\ref{sec:lag}.} We show how to obtain modules over the sheaf of symplectic functions ($\cO^\omega$ or $\cO^{\omega,\can}$) from Lagrangian branes. We need to assume the following assumption.

\medskip
\noindent\textbf{Weak unobstructedness assumption:} \textsl{For any $u\in U$, the $m_0$ term of the curved $A_\infty$ algebra on $H^*(L_u,\Lambda^\pi)$ is a scalar multiple of $\bone$, the strict unit of $H^*(L_u,\Lambda^\pi)$.} 

\medskip
\noindent This assumption implies that pairs of the form $(L_u,\alpha)$, where $L_u$ is a Lagrangian torus fiber over a point $u\in U$ and $\alpha$ is a purely imaginary one form in $H^1(L_u,\C)$, are Lagrangian branes in $\Fuk(M)$. Denote by $\Fuk^\pi(M)$ the full subcategory of $\Fuk(M)$ consisting of such Lagrangian branes. 

\begin{theorem}
\label{intro:lag}
Assume the weak unobstructedness assumption. Then there exists a linear $A_\infty$ functor $P:\Fuk^\pi(M)\ra \tw(\cO^\omega)$ which is a quasi-equivalence onto its image. The same statement also holds for $\cO^{\omega,\can}$.
\end{theorem}

\paragraph{Section~\ref{sec:koszul}.} In this section we study mirror symmetry by Koszul duality theory. We continue to work with the above assumption. First we describe the Koszul dual algebra of the sheaf $\cO^{\omega,\can}$.  For this observe that the weak unobstructedness assumption implies that for any $b\in H^1(L_u,\Lambda^\pi)$ we have
\[ m_0+m_1(b)+m_2(b,b)+m_3(b,b,b)+\cdots+m_k(b,\cdots,b)+\cdots \equiv 0 \pmod 1\]
where $1$ is the strict unit for the $A_\infty$ algebra $H^*(L_u,\Lambda^\pi)$. Then as in~\cite{FOOOtoric} we can define a potential function $W$ on the tangent bundle $TU$ by

\[ W(u, b):= \mbox{ the coefficient of $1$ of the sum\;} \sum_{k=0}^\infty m_k(b^{\otimes k}).\]
We then show that  that $W$ is in fact a ``holomorphic" function on $TU$ with its natural complex structure. Define $\cO_{TU}^\hol$ to be the ``Dolbeault complex" of the structure sheaf of $TU$ endowed with a curvature term given by $W$~\footnote{Here several words are in quote since we need to work with Novikov ring, and hence these notions need to be interpreted correctly.}. Then $\cO_{TU}^\hol$ is Koszul dual to $\cO_{M(U)}^{\omega,\can}$ in the following sense.

\begin{lemma}
\label{intro:koszul}
There is a universal Maurer-Cartan element $\tau$ in the tensor product curved $A_\infty$ algebra $ \cO_{TU}^\hol \otimes_{\Omega^*_U}\cO_{M(U)}^{\omega,\can}$.
\end{lemma}

\medskip
\noindent This universal Maurer-Cartan element can be used to construct an $A_\infty$ functor \[\Phi^\tau:\tw(\cO^{\omega,\can}) \ra \tw(\cO^\hol_{TU})\]
which is also a quasi-equivalence onto its image. Pre-composing with the functor $P$ in Theorem~\ref{intro:lag} yields the following.

\begin{theorem}
\label{intro:hms}
The composition $\Fuk^\pi(M)\stackrel{P}{\ra} \tw(\cO^{\omega,\can}) \stackrel{\Phi^\tau}{\ra} \tw(\cO^{\hol})$ is a quasi-equivalence onto its image.
\end{theorem}

\medskip
\noindent The image of this composition functor can also be described under certain assumptions.

\begin{theorem}
\label{thm:torus}
Assume that strictly negative Maslov index does not contribute to structure maps $m_k$, and assume further that the potential function $W\equiv 0$. Then the object $\Phi^\tau(\sL_u(\alpha))$ is quasi-isomorphic to a skyscraper sheaf $\Lambda^\pi(u,\alpha)$ supported at $u\in U$.
\end{theorem}

\paragraph{Section~\ref{sec:fourier}.} This section is parallel to the previous section, but replacing $\cO^{\omega,\can}$ by $\cO^{\omega}$. We also need to replace Koszul duality by its global version: Fourier-Mukai transform. Similar results as Theorem~\ref{intro:hms} and Theorem~\ref{thm:torus} are obtained. The universal Maurer-Cartan element $\tau$ in Theorem~\ref{intro:koszul} is replaced by a quantum version of the Poincar\' e bundle. Again this construction yields (local) $A_\infty$ mirror functors.  We refer to Section~\ref{sec:fourier} for more details.

\paragraph{Section~\ref{sec:functor}.} As an application of the general theory, we prove a version of homological mirror symmetry between a compact toric symplectic manifold and its Landau-Ginzburg mirror. We summarize the main results in the following theorem.

\begin{theorem}
\label{intro:toric}
Let $M$ be a compact smooth toric symplectic manifold, and denote by $\pi:M(\Delta^\interior)\ra \Delta^\interior$ the Lagrangian torus fibration over the interior of the polytope of $M$. Then there exists an $A_\infty$ functor $\Psi: \Fuk^\pi(M) \ra \tw(\cO^\hol_{T(\Delta^\interior)})$ which is a quasi-equivalence onto its image. If furthermore $M$ is Fano of complex dimension less or equal to two, then image of $\Psi$ split generates $\tw(\cO^\hol_{T(\Delta^\interior)})$ after reduction from $\Lambda^\pi$ to $\C$.
\end{theorem}

\medskip
\remark The Fano assumption is to ensure convergence over $\C$ while the dimension assumption is more of technical nature. We expect the functor $\Psi$ to be always essentially surjective over Novikov ring as long as $W$ has isolated singularities.

\paragraph{Appendices~\ref{app:modules} and~\ref{app:fm}.} We include materials on homological algebras of $A_\infty$ modules over curved $A_\infty$ algebras which might have ``internal curvatures". We also interpret Koszul duality functors as an affine version of Fourier-Mukai transform. Materials in these appendices are well-known to experts, but are not easy to find in literature. We include them here for completeness. A detailed explanation of the sign conventions used in this paper is also included in the appendix~\ref{app:fm}.

\paragraph{Acknowlegment.} The author is grateful to G. Alston, L. Amorim, K. Fukaya, and Y.-G. Oh for useful discussions on symplectic geometry; to A. Polishchuk for his help on essentially everything that involves a lattice. Thanks to D. Auroux and M. Abouzaid for sharing their comments on an earlier version of this paper. The author is also indebted to his advisor A. C\u ald\u araru for teaching him classical Fourier-Mukai transform and providing generous support during graduate school. This work is done during the author's transition from University of Wisconsin to University of Oregon, thanks to both universities for providing excellent research condition.

\paragraph{Notations and Conventions.}
\begin{itemize}
\item $\Lambda^\pi$: relative Novikov ring, see Definition~\ref{def:ring};
\item $\Fuk^\pi(M)$: the full $A_\infty$ subcategory of $\Fuk(M)$ consisting of Lagrangian torus fibers endowed with purely imaginary invariant one form;
\item $x_1,\cdots,x_n$: local action coordinates on the base $U$ of a Lagrangian torus fibration;
\item $y_1^\chk,\cdots, y_n^\chk$: angle coordinates on Lagrangian torus fibers;
\item $y_1,\cdots,y_n$: dual angle coordinates on dual torus fibers.
\item $e_1:=dy_1^\chk,\cdots, e_n:=dy_n^\chk$: the corresponding integral basis for the lattice $R^1\pi_*\Z$ over $U$;
\item We use the following two signs frequently in this paper.
\begin{itemize}
\item[(1)] $\epsilon_k:= \sum_{i=1}^{k-1} (k-i)|a_i|$ associated to homogeneous tensor products $a_1\otimes\cdots\otimes a_k$;
\item[(2)] $\eta_k=\sum_{i=1}^{k-1} |a_i|(|b_{i+1}|+\cdots+|b_k|)$ resulting from the permutation 
\[(a_1\otimes b_1)\otimes\cdots\otimes (a_k\otimes b_k) \mapsto (a_1\otimes\cdots\otimes a_k)\otimes (b_1\otimes\cdots\otimes b_k).\]
\end{itemize}
We will use two different sign conventions: one used by the authors of~\cite{FOOO} in Lagrangian Floer theory, and the other one as in~\cite{Keller}. Structure maps of the former will be denoted by $m_k$, and the latter by $m_k^\epsilon$. They are related by 
\[ m_k(a_1\otimes\cdots\otimes a_k)=(-1)^{\epsilon_k} m_k^{\epsilon}(a_1\otimes\cdots\otimes a_k).\]
\end{itemize}

\section{Symplectic functions}
\label{sec:symp}

Let $\pi:M(U)\ra U$ be a small (to be made precise below) and smooth local piece of a Lagrangian torus fibration $M\ra B$. In this section we present a construction of a sheaf of curved $A_\infty$ algebras over $U$ encoding symplectic geometry of $M$. An interesting feature of our construction is that the symplectic form $\omega$ itself enters as part of the curvature term. 

\paragraph{$A_\infty$ algebras associated to Lagrangian submanifolds.} Let $L$ be a relatively spin compact Lagrangian submanifold in a symplectic manifold $(M,\omega)$. In~\cite{FOOO} and~\cite{Fukaya} a curved $A_\infty$ algebra structure was constructed on $\Omega^*(L,\Lambda_0)$, the de Rham complex of $L$ with coefficients in certain Novikov ring $\Lambda_0$ over $\C$~\footnote{Strictly speaking this version of Floer theory was developed by K. Fukaya in~\cite{Fukaya} heavily based on the work of himself, Y.-G. Oh, H. Ohta and K. Ono~\cite{FOOO}.}. Here coefficient ring $\Lambda_0$ is defined by
\[ \Lambda_0:=\left\{\sum_{i=1}^\infty a_i T^{\lambda_i} \mid a_i\in \C, \lambda_i\in \R^{\geq 0}, \lim_{i\ra \infty} \lambda_i= \infty\right\}\]
where $T$ is a formal parameter of degree zero. Note that the map $\val:\Lambda_0\ra \R$ defined by $\val(\sum_{i=1}^\infty a_i T^{\lambda_i}):= \inf_{a_i\neq 0} \lambda_i$ endows $\Lambda_0$ with a valuation ring structure. Denote by $\Lambda_0^+$ the subset of $\Lambda_0$ consisting of elements with strictly positive valuation. This is the unique maximal ideal of $\Lambda_0$. In~\cite{FOOO} there was an additional parameter $e$ to encode Maslov index to have a $\Z$-graded $A_\infty$ structure. If we do not use this parameter, we need to work with $\Z/2\Z$-graded $A_\infty$ algebras. 

We briefly recall the construction of this $A_\infty$ structure on $\Omega^*(L,\Lambda_0)$. Let $\beta\in\pi_2(M,L)$ be a class in the relative homotopy group, and choose an almost complex structure $J$ compatible with (or simply tamed) $\omega$. Form $\cM_{k+1,\beta}(M,L;J)$, the moduli space of stable $(k+1)$-marked $J$-holomorphic disks in $M$ with boundary lying in $L_0$ of homotopy class $\beta$ with suitable regularity condition in interior and on the boundary. The moduli space $\cM_{k+1,\beta}(M,L;J)$ is of virtual dimension $d+k+\mu(\beta)-2$ (here $\mu(\beta)$ is the Maslov index of $\beta$).

There are $(k+1)$ evaluation maps $ev_i:\cM_{k+1,\beta}(M,L;J)\ra L$ for $i=0,\cdots,k$ which can be used to define a map $m_{k,\beta}:(\Omega^*(L,\C)^{\otimes k}) \ra \Omega^*(L,\C)$ of form degree $2-\mu(\beta)-k$  by formula
\[ m_{k,\beta} (\alpha_1,\cdots,\alpha_k):=(\ev_0)_!(\ev_1^*\alpha_1\wedge\cdots\wedge \ev_k^*\alpha_k).\]
To get an $A_\infty$ algebra structure we need to combine $m_{k,\beta}$ for different $\beta$'s. For this purpose we define a submonoid $G(L)$ of $\R^{\geq 0}\times 2\Z$ as the minimal one generated by the set 
\[\left\{(\int_\beta \omega,\mu(\beta))\in \R^{\geq 0}\times 2\Z \mid \beta\in \pi_2(M,L), \cM_{0,\beta}(M,L;J)\neq \emptyset\right\}.\]
Then we can define the structure maps $m_k:(\Omega(L,\Lambda_0))^{\otimes k} \ra \Omega(L,\Lambda_0)$ by
\[ m_k(\alpha_1,\cdots,\alpha_k):=\sum_{\beta\in G(L)} m_{k,\beta}(\alpha_1,\cdots,\alpha_k) T^{\int_\beta \omega}.\]
Note that we need to use the Novikov coefficients here since the above sum might not converge for a fixed value of $T$. The boundary stratas of $\cM_{k+1,\beta}(M,L;J)$ are certain fiber products of the diagram
\[\begin{CD}
\cM_{i+1,\beta_1}(M,L;J)\times_L\cM_{j+1,\beta_2}(M,L;J) @>>> \cM_{j+1,\beta_2}(M,L;J) \\
@VVV                      @VV\ev_l V\\
\cM_{i+1,\beta_1}(M,L;J) @>\ev_0>> L.
\end{CD}\]
Here $1\leq l\leq j$, $i+j=k+1$, and $\beta_1+\beta_2=\beta$. Indeed using this description of the boundary stratas the $A_\infty$ axiom for structure maps $m_k$ is an immediate consequence of Stokes formula.

We should emphasize that a mathematically rigorous realization of the above ideas involves lots of delicate constructions carried out by authors of~\cite{FOOO}. Indeed the moduli spaces $\cM_{k+1,\beta}(M,L;J)$ are not smooth manifolds, but Kuranishi orbifolds with corners, which causes trouble to define an integration theory. Even if this regularity problem is taken care of there are still transversality issues to define maps $m_{k,\beta}$ to have the expected dimension.  Moreover it is not enough to take care of each individual moduli space since the $A_\infty$ relations for $m_k$ follows from analyzing the boundary stratas in $\cM_{k+1,\beta}(M,L;J)$. Thus one needs to prove transversality of evaluation maps that are compatible for all $k$ and $\beta$. Furthermore one also need to deal with not only disk bubbles, but also sphere bubbles and regularity and transversality issues therein. We refer to the original constructions of~\cite{FOOO} and~\cite{Fukaya} for solutions of these problems. 

\paragraph{Main properties of the $A_\infty$ algebra $\Omega^*(L,\Lambda_0)$.} Let us summarize some of the main properties of $\Omega^*(L,\Lambda_0)$ proved in~\cite{FOOO} and~\cite{Fukaya}.

\begin{itemize}
\item (Invariants of symplectic geometry) The homotopy type of this $A_\infty$ structure on $\Omega(L,\Lambda_0)$ is independent of $J$, moreover it is invariant under symplectomorphism;
\item (Deformation property) This $A_\infty$ structure is a deformation of the classical differential graded algebra structure on the de Rham complex with coefficients in $\Lambda_0$~\footnote{More precisely the construction in~\cite{FOOO} yields an $A_\infty$ algebra structure whose zero energy part is only homotopic to the de Rham algebra of $L$; while in~\cite{Fukaya} the zero energy part is exactly the same.};
\item (Algebraic property) The $A_\infty$ structure can be constructed to be strict unital and cyclic so that constant function $\bone$ is the strict unit and structure maps $m_k$ are cyclic with respect to the Poincar\'e pairing $<\alpha,\beta>=\int_L \alpha\wedge\beta$.
\end{itemize}
We describe one more important property of this $A_\infty$ which is crucial for applications in this paper. This is analogous to the divisor equation in Gromov-Witten theory of closed Riemann surfaces. Such a generalization was first observed by C.-H. Cho~\cite{Cho} in the case of Fano toric manifolds. Later in~\cite{Fukaya} Cho's result was generalized by K. Fukaya to general symplectic manifolds.

\begin{lemma}[K. Fukaya~\cite{Fukaya} Lemma $13.1$ and $13.2$]
\label{lem:fukaya}
Let $b\in H^1(L,\Lambda_0)$ and consider any lift of it to an element of $\Omega^1(L,\Lambda_0)$ which we still denote by $b$. Then for any $k\geq 0$ and $l\geq 0$ we have
\[\sum_{l_0+\cdots+l_k=l} m_{k+l,\beta}(b^{\otimes l_0},\alpha_1,\cdots,\alpha_k,b^{\otimes l_k}) = \frac{1}{l!}<b,\partial\beta>^l m_{k,\beta}(\alpha_1,\cdots,\alpha_k).\]
Here $\partial\beta \in H^1(L,\Z)$ is the boundary of $\beta$.
\end{lemma}


\remark The main construction of~\cite{Fukaya} was devoted to constructing compatible Kuranishi structure and continuous multi-section perturbations that are compatible with forgetful maps between Kuranishi spaces $\cM_{k+l+1,\beta}(M,L;J)\ra \cM_{k+1,\beta}(M,L;J)$ which forget the last $l$ marked points for various $k$ and $l$. Indeed with this structure being taken care of the above lemma is a simple exercise on iterated integrals.

\paragraph{Local family of $A_\infty$ algebras.} If $\pi:M(U)\ra U$ is a local smooth family of Lagrangian torus, we get an $A_\infty$ algebra for each point $u\in U$. We would like to consider this as giving us a family of $A_\infty$ algebra over $U$. However there is a delicate point involved here: in general this fiber-wise construction does not produce an $A_\infty$ algebra over $U$ on the relative de Rham complex of $\pi$ even for a generic almost complex structure $J$ due to wall-crossing discontinuity. In~\cite{Fukaya} Section $13$ K. Fukaya constructed such a family by allowing almost complex structures to depend on the Lagrangians. Let us recall Fukaya's construction here.

Let $u\in U$ be a fixed point inside the smooth locus $B_0\subset B$. We consider a neighborhood $U$ of $u$ in $B_0$ such that there is a symplectomorphism $s: \pi^{-1}(U)\ra \cN$ identifying $\pi^{-1}(U)\subset M$ with a tubular neighborhood $\cN$ of the zero section of $T^*L_u$. Using the affine coordinates $x_1,\cdots,x_n$ and $y_1^\chk,\cdots,y_n^\chk$ to trivialize $T^*L_u$ we get an identification $\pi^{-1}(U)\cong U\times L_u$ where $U$ is considered as an open subset of $H^1(L_u,\R)=\R dy^\chk_1\oplus\cdots\oplus \R dy^\chk_n$. Equivalently this is to say that near-by Lagrangian fibers of $L_u$ maybe viewed as the graphs of one-forms on $L_u$ through the symplectomorphism $s$.

Let $p\in U$ be any point in $U$ corresponding to the one-form $p_1dy^\chk_1+\cdots+p_ndy^\chk_n$ on $L_u$. We define a \emph{diffeomorphism} $\phi_{u,p}$ of $M$ by formulas
\begin{align*}
V_{u,p} &:= p_1\partial/\partial x_1+\cdots+ p_n\partial/\partial x_n,\\
T_{u,p} &:= \epsilon_{U} \cdot (s)^{-1} V_{u,p},\\
\phi_{u,p} &:= \mbox{time one flow of $T_{u,p}$}.
\end{align*} 
Here $\epsilon_{U}$ is a cut-off function supported in a neighborhood $V$ of $U$, and is constant $1$ on $U$. The almost complex structure we shall use to form the Fukaya algebra on the fiber $L_p$ is $(\phi_{u,p})_*J$ where $J$ is the almost complex structure we use on the fiber $L_u$. By shrinking $U$ if necessary we can assume that $(\phi_{u,p})_*J$ is $\omega$-tamed for all $p\in U$.

The main advantage of this choice of almost complex structures depending on Lagrangians is that it induces identification of various moduli spaces:
\[ (\phi_{u,p})_*: \cM_{k,\beta}(M,L_u;J) \cong \cM_{k,(\phi_{u,p})_*\beta}(M,L_p;(\phi_{u,p})_*J).\]
Thus maps $m_{k,\beta}$ involved in the definition of the $A_\infty$ algebra associated a Lagrangian submanifold $L_p\; (p\in U)$ does not depend on the base parameter $p$, which implies that the structure maps $m_k:=\sum_\beta m_{k,\beta} T^{\int_\beta \omega}$ depend on the $u$-parameter only via symplectic area $\int_\beta \omega$. The follow lemma makes this dependence explicit.

\begin{lemma}
\label{lem:area}
Let $L_p$ be a near-by fiber of $L_u$ such that $L_p$ is defined as the graph of the one form $\alpha_p:=\sum_{i=1}^n p_i e_i \in H^1(L_u,\R)$. Then for each $\beta\in \pi_2(M,L_u)$ we have
\[ \int_{(\phi_{u,p})_*\beta} \omega-\int_\beta \omega = <\alpha_p,\partial \beta>\]
where $\partial\beta\in \pi_1(L_u)$ is the boundary of $\beta$.
\end{lemma}

\proof This is an exercise in Stokes' formula.

\paragraph{Relative Novikov ring.} To have a sheaf of $A_\infty$ algebras over $U$ we need to introduce another Novikov type coefficient ring. Let $\pi: M(U)\ra U$ be a small smooth local piece of Lagrangian torus fibration as above, and let $u\in U$ be the base point in $U$. The family of monoids $G(L_p) \; (p\in U)$ defines a bundle of monoids over $U$. This bundle is trivialized over $U$ to $U\times G(L_u)$ by the diffeomorphisms $\phi_{u,p}$. In other words we have a local system of monoids over $U$. We denote it by $G$.

\begin{definition}
\label{def:ring}
The relative Novikov ring $\Lambda^\pi$ associated to the family of Lagrangians $\pi:M(U)\ra U$ is defined by 
\[\Lambda^\pi:=\left\{\sum_{i=1}^\infty a_i T^{\beta_i} \mid  a_i\in \C, \beta_i\in G; \sharp\left\{a_i\mid  a_i\neq 0, \int_{\beta_i}\omega|_u\leq E\right\}<\infty, \forall E\in \R \right\}.\] 
\end{definition}

\noindent The ring $\Lambda^\pi$ is a $\Z$-graded ring with $T^\beta$ of degree $\mu(\beta)$. Note that the Maslov index map $\mu$ is well-defined on $G$ since it is preserved under isotopies.  This grading makes the operator $m_{k,\beta} T^\beta$ homogeneous of degree $2-k$.

The ring $\Lambda^\pi$ can also be endowed with a valuation by evaluating symplectic area at the base point $u\in U$, i.e.
\[\val(\sum_{i=1}^\infty a_i T^\beta):=\inf_{a_i\neq 0} { \int_{\beta_i}\omega|_u}.\]
Using this valuation map we can define a decreasing filtration on $\Lambda^\pi$ by setting $F^{\leq E}:= \val^{-1}([E,\infty))$. This filtration is called energy filtration. The ring $\Lambda^\pi$ is complete with respect to this filtration. The energy filtration is a useful tool since it induces a spectral sequence to compute Floer homology, see Chapter $6$ of~\cite{FOOO}.

\paragraph{A sheaf of $A_\infty$ algebras.} Lagrangian Floer theory developed in~\cite{FOOO} and~\cite{Fukaya} can be formulated over the ring $\Lambda^\pi$. Its relationship to the ring $\Lambda_0$ is that there is a ring homomorphism
$ \Lambda^\pi \ra \Lambda_0$ for each point $p\in U$ defined by
\[ \sum_{i=1}^\infty a_i T^{\beta_i} \mapsto \sum_{i=1}^\infty a_i T^{\int_{\beta_i}\omega}\]
where we consider $\beta_i$ as an element of $G(L_p)$.

Using this Novikov ring we can define an $A_\infty$ algebra structure on $\Omega_\pi(\Lambda^\pi)$, the relative de Rham complex with coefficients in $\Lambda^\pi$. Explicitly an element of $\Omega_\pi(\Lambda^\pi)$ is of the form $\sum_{i=1}^\infty \alpha_i T^{\beta_i}$ satisfying the same finiteness condition as in the definition of $\Lambda^\pi$. Here $\alpha_i\in \Omega_\pi(\C)$ are $\C$-valued relative differential forms. 

The structure maps of this $A_\infty$ algebra $\Omega_\pi(\Lambda^\pi)$ are defined by
\[m_k(\alpha_1,\cdots,\alpha_k):=\sum_{\beta\in G} m_{k,\beta}(\alpha_1,\cdots,\alpha_k)T^\beta,\]
and we extend these maps $\Lambda^\pi$-linearly to all elements of $\Omega_\pi(\Lambda^\pi)$. Note that since the maps $m_{k,\beta}$ does not depend on the base parameter, they are in particular smooth, which make the above definition valid.

\paragraph{D-module structure on $\Omega_\pi(\Lambda^\pi)$.} Observe that the sheaf $\Omega_{\pi}(\C)$, being the relative de Rham complex with complex coefficients, has a D-module structure over $U$.  We extend this D-module structure to $\Omega_\pi(\Lambda^\pi)$ by
\begin{equation}
\label{eq:area}
\nabla (T^\beta):= -\nabla(\int_{(F_u)_*\beta} \omega)T^\beta
\end{equation}
and Leibniz rule~\footnote{This definition is due to the fact that in the convergent case we specialize $T$ to be $e^{-1}$, as is done in~\cite{CO} Section $13$.}. We still denote this derivation by $\nabla$, and call it the Gauss-Manin connection.

\paragraph{Variational structure on $\Omega_\pi(\Lambda^\pi)$.} Our next goal is to study the variational structure of the $A_\infty$ structure $m_k$ using the Gauss-Manin connection $\nabla$. For this we consider the de Rham complex of the D-module $\Omega_\pi(\Lambda^\pi)$ over $U$. As a sheaf over $U$ this is the same as $\Omega_\pi(\Lambda^\pi)\otimes_{C^\infty_U} \Omega^*_U$. We wish to extend $m_k$ on $\Omega_\pi(\Lambda^\pi)$ to this tensor product. For this purpose it is more convenient to use a different sign convention which is better to form tensor products. We refer to the new sign convention as the $\epsilon$ sign convention. Structure maps in this sign convention will be denoted by $m_k^\epsilon$. It is related to the previous sign convention by formula
\[ m_k(a_1\otimes\cdots\otimes a_k)=(-1)^{\epsilon_k} m_k^{\epsilon}(a_1\otimes\cdots\otimes a_k)\]
where $\epsilon_k:= \sum_{i=1}^{k-1} (k-i)|a_i|$. We extend the maps $m_k^\epsilon$ on $\Omega_\pi(\Lambda^\pi)$ to its de Rham complex to get maps $m_k^\epsilon: (\Omega_\pi(\Lambda^\pi)\otimes \Omega^*_U)^k\ra \Omega_\pi(\Lambda^\pi)\otimes \Omega^*_U$ which are defined by
\begin{align*}
m_0^\epsilon&:= m_0\otimes \bone;\\
m_1^\epsilon &(f\otimes \alpha):=m_1^\epsilon(f)\otimes \alpha+(-1)^{|f|}f\otimes d_{\dR}\alpha;\\
m^\epsilon_k&((f_1\otimes \alpha_1)\otimes\cdots\otimes(f_k\otimes \alpha_k)):=(-1)^{\eta_k} m^\epsilon_k(f_1,\cdots,f_k)\otimes (\alpha_1\wedge\cdots\wedge\alpha_k)
\end{align*}
where the sign is given by $\eta_k=\sum_{i=1}^{k-1} |\alpha_i|(|f_{i+1}|+\cdots+|f_k|)$. Here we have abused the notation $m_k^\epsilon$, but no confusion should arise. It is straightforward to check that the maps $m_k^\epsilon$ defines an $A_\infty$ algebra structure on the tensor product $\Omega_\pi(\Lambda^\pi)\otimes \Omega^*_U$.

\begin{lemma}
\label{lem:diffeo}
For all $k\geq 0$ we have the following compatibility between the $A_\infty$ structure and the D-module structure on $\Omega_\pi(\Lambda^\pi)$:
\begin{equation}
\label{equ:diffeo}
[\nabla, m^\epsilon_k]=\sum_{i=1}^{k+1} (-1)^{i-1} m^\epsilon_{k+1}(\id^{i-1}\otimes \omega\otimes \id^{k-i+1}).
\end{equation}
Here $\omega$ is the symplectic form of $M$ restricted to $M(U)$, and it is viewed as an element of $\Omega_\pi(\Lambda^\pi)\otimes_{C^\infty_U} \Omega^*_U$. Locally in action-angle coordinates $\omega= \sum_{i=1}^n -dy_i\otimes dx_i^\chk=\sum_{i=1}^n -e_i\otimes dx_i^\chk$.
\end{lemma}

\proof Up to signs this lemma is a direct consequence of Lemma~\ref{lem:area} and Lemma~\ref{lem:fukaya}. We include the proof here to illustrate our sign conventions. It is enough to prove the lemma for flat sections $f_1,\cdots,f_k\in \Omega_\pi(\Lambda^\pi)$. The left hand side operator applied to $f_1\otimes\cdots\otimes f_k$ gives
\begin{align*}
&\nabla m_{k}^\epsilon (f_1\otimes\cdots\otimes f_k) = (-1)^{|f_1|+\cdots+|f_k|+2-k} \sum_\beta m_{k,\beta}^\epsilon (f_1\otimes\cdots\otimes f_k) \nabla (T^\beta)\\
&=(-1)^{|f_1|+\cdots+|f_k|+2-k} (-1)^{\epsilon_k} \sum_\beta m_{k,\beta} (f_1\otimes\cdots\otimes f_k) T^\beta \otimes \nabla(-\int_\beta \omega)\\
&= (-1)^{|f_1|+\cdots+|f_k|+1-k} (-1)^{\epsilon_k} \cdot\\
 &\cdot \sum_\beta m_{k,\beta} (f_1\otimes\cdots\otimes f_k)<\partial\beta,e_i>T^\beta\otimes dx_i^\chk \mbox{\;\;(by Lemma~\ref{lem:area})}\\
&=  \sum_{\beta, 1\leq j\leq k+1}(-1)^{|f_1|+\cdots+|f_k|+1-k} (-1)^{\epsilon_k}\cdot\\
&\cdot (-1)^{|f_1|k+\cdots+|f_{j-1}|(k-j+2)+(k-j+1)+\cdots |f_{k-1}|} \cdot m^\epsilon_{k+1,\beta} (f_1\cdots e_i \cdots f_k)T^\beta\otimes dx_i^\chk.\\
\end{align*}
The sign above can be simplified to $(-1)^{1-j+|f_1|+\cdots+|f_{j-1}|}$ which proves the lemma.

\begin{definition}[Differential $A_\infty$ algebras]
A differential $A_\infty$ algebra over a manifold $U$ is given by a triple $(E,\nabla,\omega)$ such that
\begin{itemize}
\item $E$ is a $\Z/2\Z$-graded D-module over $U$;
\item $E$ is a sheaf of $\Z/2\Z$-graded $A_\infty$ algebras over $U$;
\item $\omega$ is an even element in the de Rham complex of $E$.
\end{itemize}
Moreover these structures are compatible in the sense that equations~\ref{equ:diffeo} hold.
\end{definition}

\begin{theorem}
\label{thm:daa}
Let $(E,\nabla,\omega)$ be a differential $A_\infty$ algebra over a smooth manifold $U$. Then its de Rham complex $\Omega_U^*(E)$ also has an $A_\infty$ algebra structure. Explicitly its structure maps are given by (in the $\epsilon$ sign convention)
\begin{itemize}
\item $\hat{m}^\epsilon_0:=m^\epsilon_0-\omega$;
\item $\hat{m}^\epsilon_1:=m^\epsilon_1+\nabla$;
\item $\hat{m}^\epsilon_k:=m^\epsilon_k$ $\;\;\;$ for $k\geq 2$.
\end{itemize}
\end{theorem}

\proof This is a direct computation keeping track of the signs involved. Indeed the left hand side of equation~\ref{equ:diffeo} is the additional terms resulting from adding $\nabla$ to $m_1$ while the right hand side is exactly the terms we get by adding $\omega$ to the curvature term. The theorem is proved.

\begin{definition}[Symplectic functions]
Lemma~\ref{lem:diffeo} asserts that the triple $(\Omega_\pi(\Lambda^\pi),\nabla,\omega)$ forms a differential $A_\infty$ algebra over $U$. Theorem~\ref{thm:daa} implies that there is an $A_\infty$ algebra structure on the de Rham complex $\Omega_\pi(\Lambda^\pi)\otimes_{C^\infty_U} \Omega^*_U$. We denote this sheaf of $A_\infty$ algebras by $\cO^\omega_{M(U)}$ (or simply $\cO^\omega$), and refer to it as the sheaf of symplectic functions.
\end{definition}

\paragraph{Deformation property of $\cO^\omega$.} If $(X,\omega)$ is a symplectic manifold, then the triple $(C^\infty_X,d_{dR},\omega)$ forms a differential $A_\infty$ algebra over $X$ whose associated curved $A_\infty$ algebra is the de Rham algebra $\Omega^*_X$ endowed with a curvature term given by the symplectic form $-\omega$. This curved algebra may be thought of as ``classical" symplectic functions, and the algebra $\cO^\omega$ is a deformation of this classical algebra. The previous assertion follows from the deformation property of Fukaya algebras.

\paragraph{A variant construction.} In the end of this section we mention a variant of the sheaf $\cO^\omega$ that will be used in Section~\ref{sec:koszul}. Namely in the above constructions we could have used the canonical model which is certain minimal model $H^*(L_u,\Lambda^\pi)$ of the full de Rham complex model $\Omega^*(L_u,\Lambda^\pi)$ for each Lagrangian torus fibers. All the previous constructions go through in this case as well. 

More explicitly let $R\pi_*\Lambda^\pi$ be the push-forward of the constant sheaf $\Lambda^\pi$ via the map $\pi:M(U)\ra U$, and denote by $\ccH$ the sheaf $R\pi_*\Lambda^\pi\otimes_\C C^\infty_U$. Then $\ccH$ has a canonical Gauss-Manin connection $\nabla$ acting on it where we extend the action of $\nabla$ to $T^\beta$ by the same formula~\ref{eq:area}. Moreover the symplectic form $\omega$ can be viewed as an element in the de Rham complex of $\ccH$. Again locally in action-angle coordinates $\omega= \sum_{i=1}^n -dy_i\otimes dx_i^\chk=\sum_{i=1}^n -e_i\otimes dx_i^\chk$.
\begin{definition-theorem}
The triple $(\ccH,\nabla,\omega)$ forms a differential $A_\infty$ algebra over $U$. We shall denote by the resulting sheaf of $A_\infty$ algebras by $\cO_{M(U)}^{\omega,\can}$ (or simply $\cO^{\omega,\can}$). 
\end{definition-theorem}

\proof The proof is the same as for the triple $(\Omega_\pi(\Lambda^\pi),\nabla,\omega)$. We shall not repeat it here.

\medskip
\noindent The main advantage of $\cO^{\omega,\can}$ over $\cO^\omega$ is that the former is of finite rank over the differential graded algebra $\Omega^*_U(\Lambda^\pi)$; while its disadvantage is that its ``mirror" gives the tangent bundle of $U$ rather than the dual torus bundle, which we will discuss in Sections~\ref{sec:koszul} and~\ref{sec:fourier}. In view of results in~\cite{NZ} and~\cite{FLTZ} it is likely that $\cO^{\omega,\can}$ is related to certain equivariant symplectic geometry.

\section{From Lagrangians to modules}
\label{sec:lag}

Let $\pi:M(U)\ra U$ as in the previous section, and we continue to use notations therein. In this section we construct $A_\infty$ modules over $\cO^\omega$ from Lagrangian branes in $M$. We only consider Lagrangian branes of the form $(L_p,\alpha)$ for a Lagrangian torus fiber $L_p$  ($p\in U$) endowed with a purely imaginary torus invariant one form $\alpha$ on $L_p$. The case of general Lagrangian branes is left for future work.

Denote by $\Fuk^\pi(M)$ the full $A_\infty$ subcategory of $\Fuk(M)$ consisting of these objects. The main result of this section is the following theorem.

\begin{theorem}
\label{thm:lag}
There exists a linear $A_\infty$ functor $P:\Fuk^\pi(M)\ra \tw(\cO^\omega)$ which is a weak homotopy equivalence onto its image.
\end{theorem}

\noindent Here $\tw(\cO^\omega)$ is the $A_\infty$ category of twisted complexes over $\cO^\omega$ possibly with internal curvatures. We refer to the Appendix~\ref{app:modules} for its definition. We also remark that the theorem remains true if we replace $\cO^\omega$ by $\cO^{\omega,\can}$. This version is used in the next section to interpret mirror symmetry as Koszul duality.

\paragraph{Weak unobstructedness and potential function.} We begin to recall the notion of weak unobstructedness from~\cite{FOOO} Section $3.6$. Consider the Fukaya algebra $\Omega(L_p,\Lambda^\pi)$, and denote by $\bone$ its strict unit. An element $b\in \Omega^1(L_p,\Lambda^\pi)$ is called a weak Maurer-Cartan element if we have the equation
\[ \sum_{k=0}^\infty m_k(b^{\otimes k}) \equiv 0 \pmod \bone.\]
They are important to define Lagrangian Floer homology because they give rise to deformations of the $A_\infty$ structure on $\Omega(L_p,\Lambda^\pi)$ with square-zero differential.

In this paper we consider these elements from a more algebraic perspective: weak Maurer-Cartan elements give rise to $A_\infty$ modules with \emph{internal curvatures}. If the left hand side of the above equation were equal to zero (i.e. if $b$ is a Maurer-Cartan element), such a $b$ by definition is an $A_\infty$ module over $\Omega(L_p,\Lambda^\pi)$ (in fact this is a twisted complex). For the case of weak Maurer-Cartan elements we include relevant homological algebras in Appendix~\ref{app:modules}. What we get in this case is an $A_\infty$ module with an \emph{internal curvature}. In the following we shall freely use this notion.

On the set of weak Maurer-Cartan elements we define a function called potential function by formula
\begin{equation}
\label{eq:potential}
W(p,b):= \mbox{the coefficient of $\bone$ of the sum}\;\; \sum_{k=0}^\infty m_k(b^{\otimes k}).
\end{equation}

\begin{assumption}[Weak unobstructedness] 
\label{ass:wua}
For any $p\in U$, and any $b\in \Omega^1(L_p,\Lambda^\pi)$, the potential function
$W(p,b)$ is a scalar multiple of $\bone$, the strict unit of $\Omega(L_p,\Lambda^\pi)$.
\end{assumption}

\paragraph{Torus fibers.} Let $L_p$ be a Lagrangian torus fiber for some point $p\in U$. By the weak unobstructedness assumption~\ref{ass:wua} the pair $(L_p,0)$ is a Lagrangian brane in $\Fuk(M)$. We will construct an $A_\infty$ module $\sL_p$ over $\cO^\omega$ with an internal curvature in $C^\infty_U(\Lambda^\pi)$. Since the element $b=0$ is weakly unobstructed, it defines an $A_\infty$ module structure on $\Omega(L_p,\Lambda^\pi)$ over itself with internal curvature $W(p,0)$. The question is how to ``propagate" this structure to other points of $U$. For this we ``propagate" the weak Maurer-Cartan element $b=0$ by a differential equation using the symplectic form $\omega$. More precisely we define $\theta\in\cO^\omega$ over $U$ by the following differential equation with initial condition:
\[ \nabla \theta = \omega \mbox{\;\;\; and \;\;\;} \theta(p)=0.\]
Explicitly if $\theta$ is equal to $(x_1-p_1)e_1+\cdots+(x_n-p_n)e_n$ in local coordinates. Let us show that the element $\theta$ defines a weak Maurer-Cartan element of $\cO^\omega$. \
\begin{lemma}
\label{lem:constant}
We have $\nabla (\sum_{k=0}^\infty (-1)^{k(k-1)/2} m^\epsilon_k(\theta^k)=0$.
\end{lemma}

\proof This is a direct computation:
\begin{align*}
\nabla &(\sum_{k=0}^\infty (-1)^{k(k-1)/2} m^\epsilon_k(\theta^k) = \sum_{k=0}^\infty [\sum_{i=1}^{k+1} (-1)^{k(k-1)/2+i-1} m^\epsilon_{k+1} (\theta^{i-1},\omega,\theta^{k-i+1}) +\\
&+\sum_{j=1}^k (-1)^{k(k-1)/2+k+j-1} m^\epsilon_k (\theta^{j-1},\nabla\theta,\theta^{k-j}) ] \mbox{\;\;\; (by Lemma~\ref{lem:fukaya})} \\
&=\sum_{k=1}^\infty\sum_{i=1}^{k} (-1)^{(k-1)(k-2)/2+i-1} m_{k} (\theta^{i-1},\omega,\theta^{k-i}) - \\
&-\sum_{k=1}^\infty\sum_{j=1}^k (-1)^{k(k-1)/2+k+j-1} m_k (\theta^{j-1},\omega,\theta^{k-j})  \mbox{\;\;\; (by the equation $\nabla\theta=\omega$)} \\
&=0.
\end{align*}
For the last equality we observe that the sum $[(k-1)(k-2)/2+i-1]+[k(k-1)/2+k+i-1]$ from the signs is always odd, hence the two summations cancel out each other. The lemma is proved.

\medskip
\noindent Thus the sum $\sum_{k=0}^\infty (-1)^{k(k-1)/2} m^\epsilon_k(\theta^k)$ is a constant with respect to $\nabla$. We note that this ``constant" \emph{is not} a constant function on $U$, but rather a flat section of the D-module $\Lambda^\pi$ over $U$. Let us denote this element in $C^\infty_U(\Lambda^\pi)$ by $\cW_{(p,0)}$. By its definition we have 
\[\cW_{(p,0)}|_{p}=W(p,0).\] 
Applying the construction in Appendix~\ref{app:modules} we get a sheaf of $A_\infty$ modules over $\cO^{\omega}$ with internal curvature $\cW_{(p,0)}$. If we assume convergence of relevant power series after evaluating $T=e^{-1}$, the function $\cW_{(p,0)}$ is simply the constant function $W(p,0)$ over $U$.

\paragraph{Torus fibers with a purely imaginary closed one forms.} The above construction can also be generalized to the case $(L_p, \alpha)$ for some $\alpha\in H^1(L_p,\C)$ that is purely imaginary~\footnote{In fact we can take any element in $H^1(L_p,\Lambda_0)$. We restricted to purely imaginary ones for the purpose of mirror symmetry. See Sections~\ref{sec:koszul} and~\ref{sec:fourier}.}. This time we define $\theta\in \cO^\omega$ by
\[ \nabla\theta = \omega  \mbox{\;\;\; and \;\;\;} \theta(p)= \alpha.\]
Assumption~\ref{ass:wua}, together with the same computation as above, shows that $\theta$ defines a weak Maurer-Cartan element of $\cO^\omega$ with internal curvature given by a $\nabla$-flat section $\cW_{(p,\alpha)}\in C^\infty_U(\Lambda^\pi)$ such that $\cW_{(p,\alpha)}|_p=W(p,\alpha)$. We denote the associated $A_\infty$ module by $\sL_p(\alpha)$. 

\paragraph{Proof of Theorem~\ref{thm:lag}.} We begin to consider the endmorphism space of an object $(L_p,\alpha)$ in the Fukaya category. Recall by definition~\cite{FOOO} Section $3.6$ the Endomorphism complex $\Hom_{\Fuk(M)}((L_p,\alpha),(L_p,\alpha))$ is just the $A_\infty$ algebra $\Omega(L_p,\Lambda^\pi)$ endowed with a differential twisted by the weak Maurer-Cartan element $\alpha$. Explicitly its differential is defined as
\[ m_1^\alpha(x):=\sum_{i,j=0}^{\infty} m_{i+j+1}(\alpha^i,x,\alpha^j).\]
Similarly the complex $\Hom_{\cO^\omega}(\sL_p(\alpha),\sL_p(\alpha))$ is the $A_\infty$ algebra $\cO^\omega$ twisted by the weak Maurer-Cartan element $\theta$ associated to $(L_p,\alpha)$ as described in the previous paragraph.

For an element $\eta\in \Omega (L_p,\Lambda^\pi)$ we define an element $P(\eta)\in \cO^\omega$ by propagating $\eta$ in a flat way. Namely $P(\eta)$ is such that $\nabla(P(\eta))=0$ and $P(\eta)|_p=\eta$. We shall show that the map 
\[ P:\Hom_{\Fuk(M)}((L_p,\alpha),(L_p,\alpha))\ra \Hom_{\cO^{\omega}}(\sL_p(\alpha),\sL_p(\alpha))\]
is a linear homomorphism of $A_\infty$ algebras and a quasi-isomorphism on the underlying complexes. The fact that $\theta$ has internal curvature $\cW_{(p,\alpha)}$ implies that $P$ interchanges the curvature terms on both sides. To see that $P$ is compatible with all higher multiplications $m_k$ we need to show that
\[P[m(e^\alpha,\eta_1,e^\alpha,\cdots,e^\alpha,\eta_k,e^\alpha)]=M(e^\theta,P(\eta_1),e^\theta,\cdots,e^\theta,P(\eta_k),e^\theta)\]
where $m$ and $M$ are $A_\infty$ structures on $\Omega(L_p,\Lambda^\pi)$ and $\cO^{\omega}$ respectively, and $e^\alpha=1+\alpha+\alpha\otimes\alpha+\cdots$ is considered as an element of the bar complex of $\Omega(L_p,\Lambda^\pi)$, similarly for $e^\theta$. The two sides agree with each other at the point $p$ by definition. Hence it suffice to show that the right hand side is $\nabla$-flat. Computing $\nabla[M(e^\theta,P(\eta_1),e^\theta,\cdots,e^\theta,P(\eta_k),e^\theta)]$ using Lemma~\ref{lem:diffeo} and the condition $\nabla\theta=\omega$  shows that this is indeed the case. We omit this computation here as it is similar to the computation in Lemma~\ref{lem:constant}.

It remains to prove that $P$ is a quasi-isomorphism. In fact we will show that $P$ is a homotopy equivalence. For this we need to use the assumption that $U$ is contractible. Let $H$ be a homotopy retraction between the one-point space $p$ and $U$. It induces a deformation retraction between functions on the point $p$ (one dimensional) and the de Rham complex of $U$ with coefficients in $\C$. We denote this algebraic homotopy also by $H$. To extend $H$ to $\Omega(L_p,\Lambda^\pi)$ observe that
\[ P(\sum\eta_i T^{\beta_i})=\sum P(\eta_i)P(T^\beta_i)=\sum P(\eta_i) e^{\sum_{j=1}^d<\partial\beta_i, e_j>(x_j-p_j)} T^{\beta_i}\]
where the factor $I(\beta_i):=e^{\sum_{j=1}^d<\partial\beta_i, e_j>(x_j-p_j)}$ is by Lemma~\ref{lem:area}. Thus for each sector $T^\beta$ the operator $I(\beta)\circ H\circ I(\beta)^{-1}$ is a contracting homotopy, yielding a contracting homotopy between
\[(\Omega(L_p,\Lambda^\pi),0) \simeq(\Omega_\pi(\Lambda^\pi), \nabla).\]
Here the left hand side is endowed with a zero differential while the right hand side is endowed with only $\nabla$ as its differential. 

Next we consider the operator $M_1^\theta(-):=M(e^\theta,-,e^\theta)$ on $\Omega_\pi(\Lambda^\pi)$ as a deformation of $\nabla$ and use homological perturbation lemma to prove our theorem. If we denote by
\[ i: (\Omega(L_p,\Lambda^\pi),0) \hookrightarrow (\Omega_\pi(\Lambda^\pi), \nabla)\]
the inclusion, and 
\[ p: (\Omega_\pi(\Lambda^\pi), \nabla) \ra (\Omega(L_p,\Lambda^\pi),0)\]
the projection, then we have $pi=\id$ and $ip=\id+\nabla H+H \nabla$. By standard homological perturbation technique, the perturbed differential on $(\Omega(L_p,\Lambda^\pi))$ is given by formula
\[p M_1^\theta i + p M_1^\theta H M_1^\theta i+\cdots.\]
Observe that $M_1^\theta\circ M_1^\theta=0$ since $\theta$ is a weak Maurer-Cartan element, and that $M_1^\theta$ commutes with the homotopy $H$. These two facts imply that the induced perturbed differential on agrees with $p M_1^\theta i(-) = m_1^\alpha(-):=m(e^\alpha,-,e^\alpha)$. By the same reasoning, the perturbed homotopy 
\[ H+    H M_1^\theta H + H M_1^\theta H M_1^\theta H+\cdots\]
is the same as the old homotopy $H$, providing a homotopy equivalence between the perturbed complexes. Thus the case of endomorphisms is proved.

The case of $\Hom_{\Fuk(M)}((L_p,\alpha_1),(L_p,\alpha_2))$ is similar. Note that by its very definition in order that this $\Hom$ space is non-vanishing it is necessary to have $W(p,\alpha_1)=W(p,\alpha_2)$. Since we work over the same point $p$, the propagation map $P$ can still be defined. The previous proof carries over word-by-word to this case.

For the last case the $\Hom$ space $\Hom_{\Fuk(M)}((L_{p_1},\alpha_1),(L_{p_2},\alpha_2))$ for distinct $p_1, p_2\in U$ is zero by definition in $\Fuk(M)$. On the $\cO^\omega$-module side, the complex $\Hom_{\cO^\omega} (\sL_{p_1}(\alpha_1),\sL_{p_2}(\alpha_2))$ also has zero cohomology. This is due to the fact their internal curvatures are different, i.e. $W(p_1,\alpha_1)\neq W(p_2,\alpha_2)$. Thus the proof of Theorem~\ref{thm:lag} is finished.

\medskip
\remark\label{rmk:gq} The propagation equation $\nabla\theta=\omega$ is local in nature. Indeed when the base $B$ (or the smooth part of $B$) has nontrivial topology, there might not exist a global solution for this equation: we do not expect to be able to write the symplectic form as an exact form. 

\section{Mirror symmetry and Koszul duality}
\label{sec:koszul}
In this section we study the Koszul dual algebra of the $A_\infty$ algebra $\cO_{M(U)}^{\omega,\can}$ associated to a local Lagrangian torus fibration $\pi:M(U)\ra U$ (see Section~\ref{sec:symp} for its construction). We show that this Koszul dual algebra $\cO_{TU}^\hol$ is certain Dolbeault complex of the structure sheaf of the complex manifold $TU$ with values in the relative Novikov coefficient $\Lambda^\pi$ (possibly with a holomorphic function as curvature). Applying a well-known correspondence of modules over Koszul dual algebras gives a natural construction of an $A_\infty$ functor
\[\Phi^\tau:\tw(\cO^{\omega,\can})\ra \tw(\cO^\hol_{TU}).\]
Pre-composing with the propagation functor $P$ defined in the previous section gives a (local) mirror functor 
\[ \Phi^\tau\circ P: \Fuk^\pi(M)\ra \tw(\cO^\hol_{TU})\]
which is also an $A_\infty$ functor. We prove that this local mirror functor is a quasi-equivalence onto its images.  Under certain conditions we can also identify the images of this functor which turn out to be skyscraper sheaves of points in $TU$.

Throughout the section we continue to work with the unobstructedness assumption~\ref{ass:wua} introduced in the previous section.

\paragraph{Potential function restricted to $TU$.} Recall by equation~\ref{eq:potential} we defined a potential function on the set of weak Maurer-Cartan elements whose elements are pairs $(p,b)$ for a point $p\in U$ and an element $b\in H^1(L_p,\Lambda^\pi)$. We consider the restriction of this function to the set of pairs $(p,\alpha)$ where $\alpha$ is a purely imaginary element of $H^1(L_p,\C)$. The set of such pairs is canonically isomorphic to the tangent bundle $TU$. Locally in coordinates this identification is given by
\[(x_1,\cdots,x_n,y_1,\cdots,y_n)\mapsto (x_1,\cdots,x_d; \sum_{i=1}^n -\sqrt{-1}y_i e_i).\]
Consider the set of $\Lambda^\pi$ valued smooth functions on $TU$ which we shall denote by $C^\infty_{TU}(\Lambda^\pi)$. More precisely $C^\infty_{TU}(\Lambda^\pi)$ consists elements of the form $\sum_j f_j T^{\beta_j}$
for $f_j\in C^\infty_{TU}$, and for any energy bound $E$ there are only finitely many nonzero terms in the series. Moreover we consider $C^\infty_{TU}(\Lambda^\pi)$ as a D-module over $U$ by letting $\partial/\partial x_i$ act on it by $\db_i:=(\partial/\partial x_i+\sqrt{-1}\partial/\partial y_i)$. Recall that the operator $\partial/\partial x_i$ acts on $T^\beta$ via formula~\ref{eq:area}. Denote by $\Omega^*_U(C^\infty_{TU}(\Lambda^\pi)):=C^\infty_{TU}(\Lambda^\pi)\otimes \Omega^*_U$ the associated de Rham complex with coefficients in $C^\infty_{TU}(\Lambda^\pi)$. The differential on this de Rham complex will be denoted by $\db$.

The notation $\db$ has a geometric explanation. If we assume all convergence in Lagrangian Floer theory, we can evaluate $T$ at $e^{-1}$ and work over $\C$. Then $\Omega^*(C^\infty_{TU}(\Lambda^\pi))$ after evaluating $T$ at $e^{-1}$ can be identified with the classical Dolbeault differential $\db$ of the structure sheaf of $TU$.

\begin{lemma}
\label{lem:holo}
We have $\db W=0.$ In the geometric situation this implies $W$ is a holomorphic function on $TU$.
\end{lemma}

\proof By the identification $(y_1,\cdots,y_d)\mapsto \sum_{i=1}^n -\sqrt{-1}y_i e_i$ mentioned above, we evaluate $W$ at $b=\sum_{i=1}^n-\sqrt{-1}y_i e_i$ and differentiate. We have
\begin{align*}
\sqrt{-1}\partial_{y_i} (m_{k+1}(b,\cdots,b)) &=\sum_{j=1}^{k+1} m_{k+1}(b,\cdots,e_i,\cdots,b) \mbox{\small \;\;\; ($e_i$ in the $j$-th spot)}\\
                                                      &=-\partial_{x_i} m_k(b,\cdots,b) \mbox{\;\;\; (by Lemma~\ref{lem:diffeo}).}
\end{align*}
Summing over $k$ yields the result. Hence the lemma is proved.

\paragraph{Koszul dual of $\cO_{M(U)}^{\omega,\can}$.} The Koszul dual algebra $\cO^\hol_{TU}$ (or simply $\cO^\hol$) of $\cO^{\omega,\can}$ is defined as follows. This is a sheaf of curved differential graded algebras over $U$. Its underlying differential graded algebra is $\Omega^*_U(C^\infty_{TU}(\Lambda^\pi))=C^\infty_{TU}(\Lambda^\pi)\otimes \Omega^*_U$ endowed with the differential $\db$ and the natural tensor product algebra structure. Its curvature term is $-W$ (the sign is due to Theorem~\ref{thm:koszul} below). By the lemma above this defines a curved differential graded algebra. Geometrically we think of it as the Dolbeault complex of $TU$ endowed with a curvature term $-W$. 

The two sheaves of algebras $\cO^{\omega,\can}$ and $\cO^\hol$ are Koszul dual to each other in the sense of the following theorem.

\begin{theorem}
\label{thm:koszul}
Define the element $\tau:= \sum_{i=1}^n  -\sqrt{-1} y_i\otimes e_i  $ in an affine coordinates of $U$ (note that since $e_i$ and $y_i$ are dual to each other, the element $\tau$ is independent of coordinates). Then $\tau$ is a Maurer-Cartan element of the tensor product $A_\infty$ algebra $\cO^\hol\otimes_{\Omega_U^*(\Lambda^\pi)} \cO^{\omega,\can}$.
\end{theorem}

\proof This is a straight forward computation except that we need to pay extra attention to the signs. Indeed to form the tensor product $A_\infty$ algebra we use the $\epsilon$ sign convention as explained in more detail in the appendix. Let us denote by $M_k^\epsilon$ the resulting structure constants. Then we have
\begin{align*}
M^\epsilon_0 &= m^\epsilon_0-\omega-W;\\
M^\epsilon_1(\tau) &= m^\epsilon_1(\tau) + \nabla(\tau) + \db (\tau);\\
M^\epsilon_k(\tau) &= m^\epsilon_k(\tau,\cdots,\tau) \;\;\;\mbox{ for $k\geq 2$ .}
\end{align*}
Observe that the sum $\sum_{k=0}^\infty (-1)^{k(k-1)/2} m_k(\tau,\cdots,\tau)$ is by definition $W$. Moreover we have $\nabla(\tau)=0$ and $\db(\tau)=\sum dx_i\otimes e_i=\omega$. Thus summing over these equalities we get
\[\sum_{k=0}^\infty (-1)^{k(k-1)/2} M^\epsilon_k(\tau,\cdots,\tau)=0.\]
Thus the theorem is proved.


\paragraph{Local construction of mirror functor.} Let $A$ be an $A_\infty$ algebra and $B$ be a curved differential graded algebra over a base ring $R$. Assume both $A$ and $B$ are free $R$-modules, and that $A$ is of finite rank of $R$. Then associated to any Maurer-Cartan element $\tau\in B\otimes A$ we can form an $A_\infty$ functor $\Phi^\tau: \tw(A)\ra \tw(B)$.  Intuitively speaking the Maurer-Cartan element $\tau$ can be used to twist the tensor product $B\otimes A$ to get a rank one twisted complex over $B\otimes A$. Viewing this $B\otimes A$-module as a kernel we get a functor from $\tw(A)$ to $\tw(B)$. This construction is straightforward for ordinary algebras $A$ and $B$, but requires more explanations in the $A_\infty$ setting. Moreover we also would like to include modules with internal curvatures into this construction. Details of these homological constructions are included in the Appendix~\ref{app:fm}.

Let us apply this construction to the case $A=\cO^{\omega,\can}$ and $B=\cO^\hol_{TU}$ over the base ring $R=\Omega_U^*(\Lambda^\pi)$, which yields an $A_\infty$ functor
\[ \Phi^\tau: \tw(\cO^{\omega,\can}) \ra \tw(\cO^\hol).\]
Next we explicitly describe this local mirror functor $\Phi^\tau$. For an $A_\infty$ module $\sL$ over $\cO^{\omega,\can}$ (in particular $\sL$ must be a $\Omega_U^*(\Lambda^\pi)$-module since this is our base ring), define an $\cO^\hol$-module $\Phi^\tau(\sL)$ as follows. As a sheaf over $U$ this is just $\cO_{TU}^\hol\otimes_{\Omega_U^*(\Lambda^\pi)}\sL$. Its $\cO_{TU}^\hol$-module structure is induced from the first tensor component. On the tensor product we put a twisted differential defined by formula
\[ d:=\db\otimes\id + \sum_{k=0}^\infty \hat{\rho}_k(\tau,\cdots,\tau).\]
If we denote by $\rho_k:(\cO^{\omega,\can})^{\otimes k}\otimes_{\Omega_U^*(\Lambda^\pi)} \sL \ra \sL$ the structure constants of the $A_\infty$ module $\sL$, then the maps $\hat{\rho}_k(\tau,\cdots,\tau):\cO^\hol\otimes_{\Omega_U^*(\Lambda^\pi)}\sL\ra \cO^\hol\otimes_{\Omega_U^*(\Lambda^\pi)}\sL$ is defined by
\[\hat{\rho}_k(\tau,\cdots,\tau)(f\otimes m):= \sum_{i_1,i_2,\cdots,i_k} y_{i_1}\cdots y_{i_k}\cdot f\otimes\rho_k(e_{i_1},\cdots,e_{i_k};m).\]
The $A_\infty$ functors on $\Hom$ spaces can also be described explicitly. We refer to Appendix~\ref{app:fm} for more details.

\paragraph{Mirror of torus fibers.} Let us identify the object $\Phi^\tau(\sL_p(\alpha))$ for the $\cO^\omega$-module $\sL_p$ associated to a torus fiber $L_p$ endowed with a purely imaginary one form $\alpha\in H^1(L_p,\C)$ on it. We refer to the previous section~\ref{sec:lag} for the construction of $\sL_p(\alpha)$.

\begin{definition}
Let $p$ be a point in $U$, and let $\alpha$ be a purely imaginary one form in $H^1(L_p,\C)$. Define an $\cO^\hol_{TU}$-module $\Lambda^\pi(p,\alpha)$ as follows. As a sheaf over $U$ it is simply the skyscraper sheaf $\Lambda^\pi$ over the point $p\in U$. The $\cO^\hol_{TU}$- module structure is defined by letting an element $f\in \cO^\hol$ acting on $\Lambda^\pi$ via multiplication by $f(p,\alpha)$.
\end{definition}

\begin{proposition}
\label{prop:torus}
Assume that strictly negative Maslov index does not contribute to structure maps $m_k$ of the Fukaya algebra $H^*(L_p,\Lambda_0)$, and assume further that the potential function $W\equiv 0$ over $U$. Then the object $\Phi^\tau(\sL_p(\alpha))$ is quasi-isomorphic to $\Lambda^\pi(p,\alpha)[-n]$.
\end{proposition}

\proof Let $e_1,\cdots,e_n$ be a trivialization of $R^1\pi_*\Z$ which induces a trivialization of $R\pi_*\Z$ whose flat sections are $e_I:=e_{i_1}\wedge\cdots\wedge e_{i_j}$ for any subset $I\subset\left\{ 1,\cdots,n\right\}$. This trivialization defines an isomorphism of sheaves on $U$
 \[\Phi^\tau(\sL_p(\alpha)) \cong \prod_{I\subset\left\{1,\cdots,n\right\}}  \cO^\hol_{TU}\otimes e_I.\]
In this trivialization we can explicitly write down the differential on $\Phi^\tau(\sL_p(\alpha))$. The ``un-twisted differential" on this sheaf is simply $\db$. By Theorem~\ref{thm:duality} in Appendix~\ref{app:modules} the twisted part is given by
\[Q (f\otimes e_I):=\sum_{k\geq 0,l\geq 0} m_{k+l+1}(\tau^l,f\otimes e_I,\theta^k).\]
Here $\tau$ is as in Theorem~\ref{thm:koszul}, and recall that $\theta$ was defined in the previous section by the differential equation $\nabla\theta=\omega$ with initial condition $\theta|_p=\alpha$. Again by Theorem~\ref{thm:duality} we have
\[ (\cW_{(p,\alpha)}-W)\id +[\db,Q]-Q^2=0.\]
Since $W$ is assumed to be zero this equation simplifies to $[\db,Q]-Q^2=0$. Moreover by degree reason, since $[\db,Q]$ has Dolbeault degree one and $Q^2$ has Dolbeault degree zero,  we have $[\db,Q]=0$ and $Q^2=0$. The equation $[\db,Q]=0$ implies that $Q$ is a holomorphic operator.

To understand the cohomology of the differential $\db+Q$ on $\Phi^\tau(\sL_p(\alpha))$, we first kill the $\db$ part. For this define a morphism 
\[ F_1: (\cA^\pi\otimes \prod_{I\subset\left\{1,\cdots,n\right\}}   e_I,Q)\ra \Phi^\tau(\sL_p(\alpha))=(\cO^\hol_{TU}\otimes \prod_{I\subset\left\{1,\cdots,n\right\}}   e_I,\db+Q)\]
where $\cA^\pi$ denotes kernel of the operator $\db: C^\infty_{TU}(\Lambda^\pi) \ra C^\infty_{TU}(\Lambda^\pi)\otimes \Omega^1_U$, and the map $F_1$ is simply the embedding map $\cA^\pi\otimes e_I\hookrightarrow \cO^\hol_{TU}\otimes e_I$. The fact that $F_1$ is a quasi-isomorphism follows from the acyclicity of Dolbeault complex. We can also explicitly describe $\cA^\pi$. This is the set of of formal sums of the form
\[ \sum_\beta e^{\sum_i (<\partial \beta, e_i> x_i)} \cdot f_\beta \cdot T^\beta\]
for holomorphic functions $f_\beta$ on $TU$. These formal series satisfies the same finiteness condition as in the definition of $\Lambda^\pi$. Note that here the extra term $e^{\sum_i (<\partial \beta, e_i> x_i)}$ appears due to the non-trivial action of $\db$ on $T^\beta$ by formula~\ref{eq:area}.

The cohomology of $(\cA^\pi\otimes \prod_{I\subset\left\{1,\cdots,n\right\}}   e_I,Q)$ can be related to $\Lambda^\pi(p,\alpha)[-n]$ by defining another map
\[ F_2: (\cA^\pi\otimes \prod_{I\subset\left\{1,\cdots,d\right\}}   e_I,Q) \ra \Lambda^\pi(p,\alpha)[-n]\]
which is defined by
 \[ F_2(f \otimes e_I)=\begin{cases} 0; \mbox{\;\; if } I\neq\left\{1,\cdots,n\right\}\\
f(p,\alpha); \mbox{\;\; if } I=\left\{1,\cdots,n\right\}.\end{cases}\]
It is clear from the definition that $F_2$ is compatible with the $\cO_{TU}^\hol$-modules structures.
 
\begin{lemma}
The $\cO_{TU}^\hol$-morphism $F_2$ is a map of complexes.
\end{lemma}

\proof Below we refer to the wedge degree of $e_I$ the integer $|I|$ . Recall the structure maps $m_{k,\beta}$ are of wedge degree $2-k-\mu(\beta)$. Thus the operator $Q_{\beta}$
\[ e_I\mapsto \sum_{k\geq 0,l\geq 0} m_{k+l+1,\beta}(\tau^l, e_I,\theta^k)\]
is of wedge degree $1-\mu(\beta)$. This is because locally $\tau=\sum_{i=1}^n -\sqrt{-1}y_i e_i$ and $\theta=\sum_{i=1}^n (x_i-p_i-\sqrt{-1}\alpha_i) e_i$ which are both of wedge degree one. By our assumption negative Maslov index disks do not contribute to $m_k$, it follows that $Q_\beta$ can only increase the wedge degree when $\beta=0$ corresponding to constant maps. By constructions in~\cite{Fukaya} this zero energy part is precisely the exterior algebra on cohomology of torus~\footnote{It is essential for this proof that the zero energy part of $H^*(L_p,\Lambda^\pi)$ is given by the exterior product as opposed to an $A_\infty$ algebra homotopy equivalent to it. The latter was proved for constructions in~\cite{FOOO} Section $7.5$.}. By formula of $\tau$ and $\theta$ we get this zero energy part is
\[Q_0(e_I)=m_{2,0}(\tau,e_I)+m_{2,0}(e_I,\theta)= \sum_i (x_i+\sqrt{-1}y_i-p_i-\sqrt{-1}\alpha_i) e_i\wedge e_I\]
which vanishes at $x=p$ and $y=\alpha$. Thus we have shown that $F_2$ is a map of complexes. The lemma is proved.

\medskip
\noindent To prove the proposition it remains to show that $F_2$ is a quasi-isomorphism. For this we make use of the spectral sequence associated the energy filtration on both sides which is valid since $F_2$ preserves this filtration (we refer to~\cite{FOOO} Chapter $6$ for details of the construction of this spectral sequence). We claim that the first page of this spectral sequence is already an isomorphism. The first page is obtained as the cohomology of the operator $Q_0$. Since the factor $e^{\sum_i (<\partial \beta, e_i> x_i)}$ is non-vanishing, the cohomology for each $\beta\in G$ is the same. Thus it suffices to analyze the case when $\beta=0$. That is we would like to show that the map
\[F_0: (\prod_{I\subset\left\{1,\cdots,n\right\}}  \cA(TU) \otimes e_I, Q_0) \ra \C(p,\alpha) [-n]\]
defined by evaluating at $(p,\alpha)$ is a quasi-isomorphism where $\cA(TU)$ is the sheaf of holomorphic functions on $TU$. To this end we observe that $Q_0$ acts on the cohomology of $\db$ as the classical Koszul differential. It is well-known that the cohomology of this Koszul complex is $\C(p,\alpha)$ concentrated in degree $n$. Note that it is important here that we deal with holomorphic functions since the acyclicity of Koszul complex fails for $C^\infty$ functions. The proposition is proved.

\medskip
\remark If $M$ is a compact symplectic manifold with vanishing first Chern class (computed with a choice of almost complex structure), and $L_p$ a \emph{special} Lagrangian submanifold, then the assumption on the Maslov index is automatically satisfied. In fact in this case all the map $\mu:G(L_p)\ra 2\Z$ is identically zero.

\medskip
\remark The assumption on Maslov index is more of a technical nature while the condition that $W=0$ is necessary to describe the object $\Phi^\tau(\sL_p(\alpha))$ as a skyscraper sheaf, since otherwise $\Phi^\tau(\sL_p(\alpha))$ is only a matrix factorization of $W-W(p,\alpha)$ which is not a complex itself. In general without these assumptions we have the following theorem.

\begin{theorem}
\label{thm:hms}
The composition of $A_\infty$ functors 
\[\Phi^\tau\circ P: \Fuk^\pi(M)\ra \tw(\cO_{TU}^\hol)\] 
is a quasi-equivalence onto its image.
\end{theorem}

\proof By Theorem~\ref{thm:lag} the first functor is in fact a linear $A_\infty$ functor which is a quasi-equivalence onto its image. Thus it remains to show this for the second functor $\Phi^\tau$.  That is we would like to show that
\[\Phi^\tau: \Hom_{\cO^{\omega,\can}}(\sL_{p_1}(\alpha_1), \sL_{p_2}(\alpha_2)) \ra \Hom_{\cO^\hol}(\Phi^\tau(\sL_{p_1}(\alpha_1)),\Phi^\tau( \sL_{p_2}(\alpha_2)))\]
is a quasi-isomorphism. For this we can argue in the way as in the proof of Proposition~\ref{prop:torus}. 

In the case when $p_1\neq p_2$ or $\alpha_1\neq \alpha_2$ it was shown in the previous section that the complex $\Hom_{\cO^{\omega,\can}}(\sL_{p_1}(\alpha_1), \sL_{p_2}(\alpha_2))$ has zero cohomology. On the other hand the complex $\Hom_{\cO^\hol}(\Phi^\tau(\sL_{p_1}(\alpha_1)),\Phi^\tau( \sL_{p_2}(\alpha_2)))$ also has vanishing cohomology. For this we can use the spectral sequence associated to the energy filtration to calculate its cohomology. The first page of this spectral sequence is already zero due to the vanishing of $\Ext^*_{TU}(\cO_{(p_1,\alpha_1)},\cO_{(p_2,\alpha_2)})$.

Thus in the following we consider the case when $p_1=p_2=p$ and $\alpha_1=\alpha_2=\alpha$. We can kill the $\Omega_U^*$ part by observing the following commutative diagram.
\[\begin{CD}
\Hom_{\Fuk^\pi(M)} ((L_p,\alpha),(L_p,\alpha)) @>\phi^\tau>> \cA^\pi\otimes \End_\C(H^*(L_p,\C))\\
@V P VV                @VVV \\
\Hom_{\cO^{\omega,\can}}(\sL_{u}(\alpha), \sL_{u}(\alpha)) @>\Phi^\tau>> \Hom_{\cO^\hol}(\Phi^\tau(\sL_{u}(\alpha)),\Phi^\tau( \sL_{u}(\alpha)))
\end{CD}\]
Several explanations of this diagram are in order. First of all as in the proof of Proposition~\ref{prop:torus}, $\cA^\pi$ is holomorphic functions on $TU$ with values in $\Lambda^\pi$, and the right vertical map is the inclusion 
\[ \cA^\pi\otimes \End_\C(H^*(L_p,\C))\hookrightarrow \cO^\hol\otimes \End_\C(H^*(L_p,\C))\]
which is a quasi-isomorphism proved in Proposition~\ref{prop:torus}. 

Secondly the left vertical arrow is the propagation map defined in Theorem~\ref{thm:lag} which is also a quasi-isomorphism. Thus in order to show the bottom arrow $\Phi^\tau$ is a quasi-isomorphism it suffices to define the top arrow $\phi^\tau$ and show that it is also a quasi-isomorphism.

We define the map $\phi^\tau$by
\[ \phi^\tau(e_I):= 1\otimes e_J^\chk \otimes \sum_{l\geq 0,i_0\geq 0,i_1\geq 0} m_{l+2+i_0+i_1}(\tau^l,e_J,\theta^{i_0},e_I,\theta^{i_1})\]
where as before $\tau=\sum -\sqrt{-1}y_i e_i$ and $\theta=\sum (x_i-p_i-\sqrt{-1}\alpha) e_i$. This formula is the same as the formula to define $\Phi^\tau$, and hence the above diagram is commutative. See Appendix~\ref{app:fm} for the formula for $\Phi^\tau$ on $\Hom$-spaces.

It remains to show that $\phi^\tau$ is a quasi-isomorphism. For this we consider the spectral sequences associated energy filtrations on both sides. The first page of this spectral sequence is already an isomorphism, again it suffices to prove that the zero energy part sector (corresponding to $\beta=0$) is an isomorphism, which reduces to show that the map of complexes
\begin{align*}
 (H^*(L_p,\C),0)&\ra (\End_{\cA(TU)}(\cA(TU)\otimes H^*(L_p,\C)), [Q_0,-])\\
e_I & \mapsto (f\otimes e_J \mapsto f\otimes e_I\wedge e_J)
\end{align*}
where $\cA(TU)$ is holomorphic functions on $TU$ and $Q_0$ is the Koszul differential associated to the regular sequence $\left\{x_i+\sqrt{-1}y_i-p_i-\sqrt{-1}\alpha_i\right\}_{i=1}^n$, is a quasi-isomorphism. This an exercise in classical Koszul duality theory between exterior algebras and symmetric algebras. The theorem is proved.

\medskip
\remark Without the assumptions in Proposition~\ref{prop:torus} it is not clear how to identify the image of $\Phi^\tau$ in $\tw(\cO^\hol)$. In general these objects are matrix factorizations of $W$ minus a constant (the internal curvatures), and we expect them to be homotopy equivalent to stabilization matrix factorizations introduced in~\cite{Dyck} Section $2.3$. The trouble to prove this claim is that the matrix factorization defined by the operator $Q (f\otimes e_I):=\sum_{k\geq 0,l\geq 0} m_{k+l+1}(\tau^l,f\otimes e_I,\theta^k)$ mixes various the wedge degrees of $e_I$ while the stabilizations introduced in~\cite{Dyck} only involves operators of wedge degrees $1$ and $-1$.

\section{Mirror symmetry and Fourier-Mukai transform}
\label{sec:fourier}

From the symplectic point of view the two sheaves of $A_\infty$ algebras $\cO^{\omega,\can}$ and $\cO^\omega$ are not much different since $\cO^{\omega,\can}$ in some sense is the minimal model of $\cO^\omega$. However from the point of view of mirror symmetry there is an important difference between the two. Indeed we have seen from the previous section that the mirror of $\cO^{\omega,\can}$ is $TU$ while in this section we show the mirror of $\cO^\omega$ is the dual torus bundle $M^\chk(U)$. Analogous results such as Proposition~\ref{prop:torus} and Theorem~\ref{thm:hms} are obtained in this case as well.

Throughout the section we continue to work with the unobstructedness assumption~\ref{ass:wua}.

\paragraph{The case without quantum corrections.} Let us first consider the case when there exists no non-trivial holomorphic disks in $M$ with boundary in $L_u$ for all $u\in U$. In this case we can work over $\C$. We shall see how relative (over U) Poincar\'e bundles appear here.

In Theorem~\ref{thm:koszul} we constructed a Maurer-Cartan element in the tensor product algebra $ \cO_{TU}^\hol\otimes \cO^{\omega,\can}$ which is of the form $\tau:=\sum_{i=1}^n -\sqrt{-1}y_i\otimes e_i$. Considering elements $e_i$ as translation invariant one forms in $\Omega^1(L_u)$, the element $\tau$ can be viewed as an element in $\cO_{TU}^\hol\otimes\cO^\omega$. As is explained in Appendix~\ref{app:fm} Koszul duality can be considered as a special case of Fourier-Mukai transform. Indeed the kernel to construct the Koszul functor $\Phi^\tau:\tw(\cO^\omega)\ra \Tw(\cO_{TU}^\hol)$~\footnote{Here we need to use the capital $\Tw$ since $\cO^\omega$ is of infinite rank over $\Omega_U^*(\Lambda^\pi)$.} is simply the rank one twisted complex over $\cO_{TU}^\hol\otimes\cO^\omega$ defined by the Maurer-Cartan element $\tau$. 

We can explicitly describe this twisted complex in coordinates $x$, $y$ and $y^\chk$. Recall $x$ is coordinates on the base $U$; $y^\chk$ and $y$ are angle coordinates on Lagrangian torus and dual torus respectively. We also trivialize $M(U)$ to $T\times U$, and denote by $H_1(T,\R)$ by $V$, its dual by $V^\chk$. Then the tensor product $\cO_{TU}^\hol\otimes\cO^\omega$ is simply $C^\infty_{V^\chk}\otimes \Omega^*(T) \otimes \Omega^*_U$~\footnote{Here and in the following $\otimes$ means completed tensor product.} endowed with the differential $d_T+\db$. Here $d_T$ is the de Rham differential on $\Omega^*(T)$. The operator $\db$ is defined in the previous section as $\sum_{i=1}^n(\partial/\partial x_i+\sqrt{-1}\partial/\partial y_i) dx_i$. Using the Maurer-Cartan element $\tau$ this differential is twisted to $d_T-\sqrt{-1}y_i e_i +\db$ where the part additional operator stands for multiplication by $-\sqrt{-1}y_ie_i$. Note that the square of this twisted operator is not zero but the symplectic form $\omega$. Indeed this should be a twisted complex over $\cO_{TU}^\hol\otimes\cO^\omega$ which has curvature $-\omega$. We shall denote this twisted complex by $K^\tau$.

To pass from $TU=V^\chk\times U$ to $M^\chk(U)=T^\chk\times U$ it suffices to quotient out the dual lattice group $\Gamma:=(H_1(T,\Z))^\chk=H_1(T^\chk,\Z)$ in $V^\chk$. However we note that the kernel does not descend to this quotient in an obvious way since the twisted operator $d_T-\sqrt{-1}y_ie_i +\db$ is not $\Gamma$-equivariant under the natural translation action. This is where Poincar\'e bundle comes into play: we can define another $\Gamma$ action on $K^\tau$ so that the operator $d_T-\sqrt{-1}y_ie_i +\db$ becomes equivariant. This ``twisted action" is given by
\begin{equation}
\label{eq:action}
 \gamma[f(y^\chk,y)]:= e^{\sqrt{-1}\gamma\cdot y^\chk} f(y^\chk,y-\gamma)
\end{equation}
where $\gamma\cdot y$ is the natural pairing between $V^\chk$ and $V$. It is well-known that if we take the above action and consider invariants in the function part $C^\infty_{V^\chk}\otimes C^\infty_T \otimes C^\infty_U$ of $K^\tau$ we get $C^\infty$-sections of the relative Poincar\'e bundle $\cP$ on $M(U)\times_U M^\chk(U)$. A direct computation verifies the following commutative diagram
\[\begin{CD}
K^\tau @>d_T-\sqrt{-1}y_ie_i +\db>>  K^\tau \\
@V\gamma VV                             @VV\gamma V \\
K^\tau @>d_T-\sqrt{-1}y_ie_i +\db>> K^\tau.
\end{CD}\]
Thus the operator $d_T-\sqrt{-1}y_ie_i +\db$ descends to an operator on invariants $(K^\tau)^\Gamma$. Moreover the $\cO^\hol_{TU}$-module structure
\[ \cO^\hol_{TU}\otimes K^\tau \ra K^\tau\]
is $\Gamma$-equivariant if we put the ordinary translation action on $\cO^\hol_{TU}$ and the twisted action on $K^\tau$. Taking invariants yields an action
\[ \cO^\hol_{M^\chk(U)} \otimes (K^\tau)^\Gamma \ra (K^\tau)^\Gamma.\]
where $\cO_{M^\chk(U)}^\hol$ is the $\Gamma$-invariants of $\cO^\hol_{TU}$ under translation action, i.e. it is the Dolbeault complex of the structure sheaf of the complex manifold $M^\chk(U)$. We can then consider $(K^\tau)^\Gamma$ as a ``twisted complex" over $\cO^\hol_{M^\chk(U)}\otimes \cO^\omega$. Here we put twisted complex in quote since the two objects $(K^\tau)^\Gamma$ and $\cO^\hol_{M^\chk(U)}\otimes \cO^\omega$ are not isomorphic even as sheaves on $U$~\footnote{Indeed the Poincar\'e bundle $\cP$ is not a topologically trivial vector bundle.}. As in the previous section we would like to use $(K^\tau)^\Gamma$ as a Kernel to define a local mirror symmetry functor. We shall perform this construction in the general case when quantum corrections are presented.

\paragraph{Fourier transform for families.} We briefly recall Fourier transform for families which may be thought of as mirror symmetry in the case when torus fibers do not bound any nontrivial holomorphic disk. Our main references are articles~\cite{AP},\cite{LYZ} and~\cite{BMP}. 

For a real torus $T$, and its dual torus $T^\chk$, there exists a Poincar\' e bundle $\cP$ over $T\times T^\chk$. The bundle $\cP$ has a natural connection $\nabla$ whose curvature is equal to the canonical symplectic form on $T\times T^\chk$. If we split the connection $\nabla$ into $\nabla^{1,0} +\nabla^{0,1}$ corresponding to the splitting $\Omega_{T\times T^\chk}^1=p_1^*\Omega_T^1\otimes p_2^*\Omega_{T^\chk}^1$, then the partial connection $\nabla^{1,0}$ has zero curvature. Thus we may form the partial de Rham complex $p_1^*\Omega_T\otimes \cP$ of the connection $\nabla^{1,0}$. Using $p_1^*\Omega_T\otimes \cP$ as the kernel of Fourier-Mukai transform, one can show that there is an equivalence between the category of local systems on $T$ and the category of skyscraper sheaves on $T^\chk$. See for example Propositions $2.6$, $2.7$ and $2.8$ in~\cite{BMP}.

In the relative case, one considers a family of Lagrangian torus $X\ra B$, and its dual family $X^\chk\ra B$. Again we main form $\cP$ over $X\times_B X^\chk$, and a partial flat connection $\nabla^{1,0}$. Using the partial de Rham complex of $(\cP,\nabla^{1,0})$ as the integral kernel, one deduces an equivalence between the category of local systems supported on fibers of $X\ra B$ and the category of skyscraper sheaves on $X^\chk$. See Section $3.2$ in~\cite{BMP}.

\paragraph{Quantum Fourier-Mukai transform.} Let us return to the general case to allow non-trivial holomorphic disks to enter the picture. The previous discussion motivates us to perform the construction in two steps:
\begin{itemize}
\item[\Rnum{1}.] Replace the sheaf $\cO^{\omega,\can}$ by $\cO^\omega$, and perform the same construction as in the previous section;
\item[\Rnum{2}.] Descend to $\Gamma$-invariants by action~\ref{eq:action}.
\end{itemize}
We begin with Step $\Rnum{1}$ which is almost word by word as in the previous section. Consider the $\cO^\omega$-module $\sL_p(\alpha)$ associated to a Lagrangian $L_p$ endowed with a purely imaginary closed one form $\alpha$ (see Section~\ref{sec:lag} for the construction of $\cO^\omega$). This is a twisted complex of rank one over $\cO^\omega$ defined by a Maurer-Cartan element $\theta$ such that $\nabla\theta=\omega$ and $\theta|_p=\alpha$. In local coordinates $\theta=\sum_i (x_i-p_i+\alpha) e_i$, and let $\tau:=\sum_i -\sqrt{-1} y_i e_i \in \cO_{TU}^\hol\otimes \cO^\omega$ be as before. By constructions in the previous section we get a twisted complex $\Phi^\tau(\sL_p(\alpha))$ over $\cO_{TU}^\hol$ which, as a $\cO_{TU}^\hol$-module, is simply $\cO_{TU}^\hol\otimes \cO^\omega$. It is endowed with a twisted differential $d:=\db+Q$ where the operator $Q$ is
\[ Q(?)=\sum_{k\geq 0,l\geq 0} m_{k+l+1}(\tau^l,?,\theta^k)\]
expressed using structure maps on the tensor product $A_\infty$ algebra $\cO_{TU}^\hol\otimes \cO^\omega$.

For the step $\Rnum{2}$ we first observe that it follows from Lemma~\ref{lem:fukaya} the potential function $W$, \emph{A priori} defined on $TU$, is in fact a function on $M^\chk(U)$ with values in $\Lambda^\pi$. This enables us to define the sheaf of curved algebras $\cO^\hol_{M^\chk(U)}$ in general. Our next goal is to show that the operator $d:=\db+Q$ on $\cO_{TU}^\hol\otimes \cO^\omega$ intertwines with action~\ref{eq:action}.

For this, we shall choose an almost complex structure $J$ on $M$ so that when restricted to $\pi^{-1}(U)$ it is given by a standard one in action-angle coordinates. Namely, over $\pi^{-1}(U)$, in action-angle coordinates $x_1,\cdots, x_n$, $y_1^\chk,\cdots, y_n^\chk$, the almost complex structure is
\[ J: \partial/\partial x_i \mapsto \partial /\partial y_n^\chk.\]
The reason we need this is because the boundary of a $J$-holomorphic disk is geodesic in the metric $\omega(-,J-)$, and with this choice of $J$, geodesics are straight lines. Thus if we consider the evaluation maps
\begin{align*}
\ev_1& : \cM_{1+1,\beta} \ra L,\\
\ev_0&: \cM_{0+1,\beta}\ra L,
\end{align*}
and the forgetful map
\[ \cM_{1+1,\beta}\ra \cM_{0+1,\beta},\]
we have the identity
\begin{equation}~\label{geodesic:eq}
\ev_1([u,t])=\ev_0([u])+t\cdot \partial\beta
\end{equation}
where $t$ is coordinate on the fiber of the forgetful map. We use this identification in the proof of the following lemma.

\begin{lemma}
\label{lem:intertwine}
For each $\beta\in G$, and $f\in \cO_{TU}^\hol\otimes \cO^\omega$ we have
\[\sum_{l\geq 0} m_{l+1+k,\beta}(\tau^l,\gamma(f),\theta^k)= \gamma [\sum_{l\geq 0} m_{l+1+k,\beta}(\tau^l,f,\theta^k)].\]
\end{lemma}

\proof We denote by $T_\gamma$ the translation action $y\mapsto y-\gamma$. Thus $\gamma(f)=e^{\sqrt{-1}\gamma\cdot y^\chk}T_\gamma^*(f)=e^{\sqrt{-1}\gamma\cdot y^\chk} f(y^\chk,y-\gamma)$. We have
\begin{align*}
& m_{l+1+k,\beta}  (\tau^l,\gamma(f),\theta^k) = (\ev_0)_! [\int_{0\leq t_1\leq t_2\cdots\leq t_{k+l+1}\leq 1}  \ev^*(\tau^l) \cdot \ev_{l+1}^*(\gamma(f))\cdot \ev^*(\theta^k)] \\
&=(\ev_0)_! [ \frac{1}{l!} <\partial\beta,\tau>^l\int_{0\leq t_{l+1}\cdots\leq t_{k+l+1}\leq 1} t_{l+1}^l\ev_{l+1}^*(\gamma(f))\cdot \ev^*(\theta^k)]\\
&=(\ev_0)_! [\frac{1}{l!}<\partial\beta,\tau>^l \int_{0\leq t_1\cdots\leq t_{k+1}\leq 1}
t_1^l\ev_1^*(e^{\sqrt{-1}\gamma\cdot y^\chk})\ev_1^*(T_\gamma f)\cdot\ev^*(\theta^k)]\\
&=(\ev_0)_! [\ev_0^*(e^{\sqrt{-1}\gamma\cdot y^\chk}) \frac{<\partial\beta,\tau>^l}{l!} \int_{0\leq t_1\cdots\leq t_{k+1}\leq 1}
t_1^l e^{\sqrt{-1}\gamma\cdot (t_1\partial\beta)}\ev_1^*(T_\gamma f)\cdot\ev^*(\theta^k)].
\end{align*}
The last equality follows from Equation~\ref{geodesic:eq}. Apply projection formula $f_!(f^*a\cdot b)=a \cdot f_! b$ to the last summation yields
\[e^{\sqrt{-1}\gamma\cdot y^\chk} (\ev_0)_! [\frac{1}{l!}<\partial\beta,\tau>^l \int_{0\leq t_1\cdots\leq t_{k+1}\leq 1}
t_1^l e^{\sqrt{-1}\gamma\cdot (t_1\partial\beta)}\ev_1^*(T_\gamma f)\cdot\ev^*(\theta^k)].\]
Summing over $l$ yields
\begin{align*}
&\sum_{l\geq 0} m_{l+1+k,\beta}  (\tau^l,\gamma(f),\theta^k)=\\
&= e^{\sqrt{-1}\gamma\cdot y^\chk} (\ev_0)_! [\int_{0\leq t_1\cdots\leq t_{k+1}\leq 1}
e^{(t_1\partial\beta)\cdot \tau} e^{\sqrt{-1}\gamma\cdot (t_1\partial\beta)}\ev_1^*(T_\gamma f)\cdot\ev^*(\theta^k)]\\
&= e^{\sqrt{-1}\gamma\cdot y^\chk} (\ev_0)_! [\int_{0\leq t_1\cdots\leq t_{k+1}\leq 1}
e^{(t_1\partial\beta)\cdot T_\gamma\tau}\ev_1^*(T_\gamma f)\cdot\ev^*(\theta^k)]\\
&= \gamma\left\{ (\ev_0)_! [\int_{0\leq t_1\cdots\leq t_{k+1}\leq 1}
e^{(t_1\partial\beta)\cdot\tau}\ev_1^*(f)\cdot\ev^*(\theta^k)]\right\}\\
&=\gamma\left\{ \sum_{l\geq 0} (\ev_0)_! [\int_{0\leq t_1\cdots\leq t_{k+1}\leq 1}
\frac{1}{l!} <\partial\beta,\tau>^l t_1^l\cdot \ev_1^*(f)\cdot\ev^*(\theta^k)]\right\}\\
&=\gamma [\sum_{l\geq 0} m_{l+1+k,\beta}(\tau^l,f,\theta^k)].
\end{align*}
Thus the lemma is proved.

\medskip
\remark It follows from the proof that in the above formula the part $\theta^k$ can be replaced by any elements. Also, the proof we gave above assumes transversality for moduli spaces involved. At present there are several different approaches to deal with transversality issues. In any case, we believe that geometric argument as in the proof of the lemma should work in full generality with any successful approach towards solving the transversality problems.

\medskip
\noindent Returning to the discussion of the twisted complex $\Phi^\tau(\sL_p(\alpha))$ endowed with differential $\db+Q$, the above lemma implies that the operator $Q$ is $\Gamma$-equivariant. Since $\db$ is also equivariant the sum operator $\db+Q$ is also $\Gamma$-equivariant. So we can take $\Gamma$-invariants to get a $\cO^\hol_{M^\chk(U)}$-module structure on the complex $[\Phi^\tau(\sL_p(\alpha))]^\Gamma$. The generalization of this construction to general objects of $\tw(\cO^\omega)$ is straight-forward. Thus we have describe a functor from $\tw(\cO^\omega)$ to $\Tw(\cO^\hol_{M^\chk(U)})$ on the level of objects. The capital $\Tw$ stands for twisted complexes of possibly infinite rank. This is necessary for us here. We denote this functor by $\Phi^\cP$ where we think of $\cP:=[K^\tau]^\Gamma$ as a certain quantized relative Poincar\'e bundle. 

We continue to define $\Phi^\cP$ on morphisms. In fact the functor $\Phi^\cP$ can be constructed as an $A_\infty$ functor $\tw(\cO^\omega)\ra \Tw(\cO^\hol_{M^\chk(U)})$. Indeed by Appendix~\ref{app:fm} the $A_\infty$ functor $\tw(\cO^\omega) \ra \Tw(\cO^\hol_{TU})$ associated to the twisting cochain $\tau$ has the form
\[\Phi^\tau(a_1,\cdots,a_k)(x):= \sum_{l\geq 0,i_0\geq 0,\cdots,i_k\geq 0} m_{l+k+1+i_0+\cdots+i_k}(\tau^l,x,\theta^{i_0},a_1,\theta^{i_1},\cdots,a_k,\theta^{i_k})\]
which by Lemma~\ref{lem:intertwine} (and its following remark) is also $\Gamma$-equivariant. Hence this functor descends to $\Gamma$-invariants to give the desired $A_\infty$ functor $\Phi^\cP$.

\paragraph{Mirror dual of torus fibers.} Since the image of $\Phi^\cP$ are twisted complexes of infinite rank over $\cO^\hol_{M^\chk(U)}$, it is \emph{A priori} unclear whether these objects are quasi-isomorphic to an object in the bounded derived category of coherent sheaves on $M^\chk(U)$ (with coefficients in $\Lambda^\pi$). We show this is the case for Lagrangian torus fibers endowed with a unitary line bundle $[\alpha]$ on it. Here the $\alpha$ is a connection one form corresponding to the line bundle $[\alpha]$ which is determined up to translation by lattice points.

\begin{proposition}
\label{prop:image}
Assume that strictly negative Maslov index does not contribute to structure maps $m_k$, and assume further that the potential function $W\equiv0$ over $U$. Then the object $\Phi^\cP(\sL_p(\alpha))$ is quasi-isomorphic to the skyscraper sheaf $\Lambda^\pi(p,\alpha)[-n]$ over the point $p\in U$. As in Proposition~\ref{prop:torus} we let an element $f\in \cO_{M^\chk(U)}^\hol$ act on $\Lambda^\pi$ via multiplication by $f(p,\alpha)$.
\end{proposition}

\proof The proof is analogous to that of Proposition~\ref{prop:torus}. Let us trivialize $M(U)=T\times U$ which induces an identification 
 \[\Phi^\cP(\sL_p(\alpha))\cong [(C^\infty_{V^\chk}\otimes \Omega^*_U) \otimes \Omega^*(T)]^\Gamma.\]
This complex is endowed with the $\Gamma$-equivariant differential $\db+Q$ as describe above. Moreover we have $[\db,Q]=Q^2=0$ as before. To calculate the cohomology of this complex we first observe a quasi-isomorphism
\[ F_1: ( [\cA(TU)\otimes \Omega^*(T)]^\Gamma, Q) \ra ([(C^\infty_{V^\chk}\otimes \Omega^*_U) \otimes \Omega^*(T)]^\Gamma,\db+Q)\]
where recall that $\cA(TU)$ is $\Lambda^\pi$ valued holomorphic function on $TU$. The map $F_1$ is defined by 
\[ [(f\cdot T^\beta )\otimes \zeta] \mapsto (e^{\sum_i (<\partial \beta, e_i> x_i)} \cdot f \cdot T^\beta)\otimes \zeta\]
where again the extra term $e^{\sum_i (<\partial \beta, e_i> x_i)}$ appears due to the non-trivial action of $\db$ on $T^\beta$~\ref{eq:area}. That $F_1$ is a quasi-isomorphism follows from the exactness of Dolbeault complex.

Next we define a morphism of $\cO_{M^\chk(U)}^\hol$-modules
\[ F_2: ( [\cA(TU)\otimes \Omega^*(T)]^\Gamma, Q) \ra \Lambda^\pi(p,\alpha)[-n] \]
by formula $F(f\otimes \zeta):=f(p,\alpha)\int_T \zeta$. It is clear that $F_2$ respects the action of $\cO_{M^\chk(U)}^\hol$ on both sides. Let us check that $F_2$ is a morphism of complexes, i.e. we would like to show that $F\circ Q=0$. For this observe that integration on $T$ kills all elements of form degree strictly less than the dimension $n$ of $T$. Moreover by the Maslov index assumption the only operator that increases this degree is $Q_0$ corresponding to trivial holomorphic disks. As in Proposition~\ref{prop:torus} this operator is explicitly given by
\begin{align*}
Q_0(f\otimes \zeta) &=f\otimes d_{\dR} \zeta+ m_{2,0}(\tau,f\otimes \zeta)+m_{2,0}(f\otimes\zeta,\theta)\\
&= f\otimes d_{\dR} \zeta + \sum_i (x_i+\sqrt{-1}y_i-p_i-\sqrt{-1}\alpha_i)\cdot f \otimes e_i\wedge \zeta.
\end{align*}
Applying $F_2$ to this sum, which by definition is evaluation at $(p,\alpha)$ and integrate over $T$, yields zero. Thus we have shown that $F_2$ is a map of complexes. It remains to prove that it is also a quasi-isomorphism. For this we consider the spectral sequences associated to energy filtration on both sides. As in the proof of Proposition~\ref{prop:torus} it suffices to analyze the case for $\beta=0$. This is done in the following lemma, which finishes the proof the proposition.

\begin{lemma}
Denote by $\cA(TU,\C)$ the sheaf of holomorphic functions on $TU$ with values in $\C$, then the map
\[ \phi: ([\cA(TU,\C) \otimes \Omega^*(T)]^\Gamma,Q_0) \ra \C[-n]\]
defined by $\phi(f\otimes \zeta):=f(p,\alpha)\cdot \int_T \zeta$ is a quasi-isomorphism.
\end{lemma}

\proof This is a classical result in Fourier transform for families. We include a proof here for completeness. The proof is similar to that of Proposition $2.6$, $2.7$ and $2.8$ in~\cite{BMP}. We only do this for the case when the dimension of $T$ is one, i.e. $T$ is a circle. The general case follows from Kunneth type argument.

We work in the universal cover of $T$ with affine coordinate $y^\chk$. Coordinates on $TU$ are $x$ and $y$. The operator $Q_0$ acts by $[\partial/\partial y^\chk + (z-p-\sqrt{-1}\alpha)]dy^\chk$ where $z=x+\sqrt{-1}y$. To analyze the cohomology of $Q_0$ it is better to work in another ``gauge", i.e. we conjugate the operator $Q_0$ with an automorphism which is given by multiplication by $e^{(z-p-\sqrt{-1}\alpha)\cdot y^\chk}$. Under this conjugation the operator $Q_0$ is identified with $\partial/\partial y^\chk \circ dy^\chk$, the de Rham differential in $y^\chk$-direction. Moreover the conjugation also changes the lattice group action on the variables $y^\chk$ and $y$. In the $y$-direction $\Gamma\subset \R^\chk$ acts simply by translation, while in the $y^\chk$-direction $m\in \Z\subset \R$ acts on $s\in C^\infty(\R\times \R^\chk)$ by
\[ s(y^\chk,y) \mapsto e^{m\cdot(z-p-\sqrt{-1}\alpha)} s(y^\chk-m,y).\]
To prove the lemma it suffices to show that an element $f(z)\otimes s(y^\chk) dy^\chk$ is exact if and only if $f(p,\alpha)\int_0^1 s(v)dv=0$. The element $f\otimes s dy^\chk$ is exact if there exists an anti-derivative of the form
\[   f(z)\otimes t(y^\chk):= f(z) \otimes \int_0^{y^\chk}  s(v) dv + h(z)\]
which is periodic in $y$-direction and in $y^\chk$-direction we have
\[ f(z) e^{m\cdot(z-p-\sqrt{-1}\alpha)}\otimes t(y^\chk-m)= f(z)\otimes t (y^\chk).\]
Using the fact that $f\otimes s$ is lattice group invariant we find that such an anti-derivative exists if and only if
\[ (e^{m\cdot(z-p-\sqrt{-1}\alpha)}-1) h= f\int_0^ms(v)dv\]
for all $m$. Denote by $q:=e^{z-p-\sqrt{-1}\alpha}$, then using the lattice group invariance of $s(v)$ we get $\int_0^m s(v)dv=(1+q+\cdots+q^{m-1})\int_0^1 s(v)dv$. 

If $f(p,\alpha)\int_0^1 s(v)dv=0$, then either $f(p,\alpha)=0$ or $\int_0^1 s(v) dv=0$. In the first case the fraction
\[ \frac{f\int_0^ms(v)dv}{q^m-1}= \frac{f\int_0^1 s(v)dv}{q-1}\]
which is independent of $m$ extends to the point $(p,\alpha)$ since $(e^{m\cdot(z-p-\sqrt{-1}\alpha)}-1)$ vanishes in first order at $(p,\alpha)$. If $\int_0^1s(v)dv=0$, then $\int_0^m s(v)dv=0$. Hence we can take $h$ to be simply zero. So in either case the form $f\otimes s dy^\chk$ is exact.

Conversely if $f\otimes s dy^\chk$ is exact then $f(p,\alpha)\int_0^m s(v)dv=0$ of all $m$. In particular we have $f(p,\alpha)\int_0^1s(v)dv=0$. The lemma is proved.

\medskip
\noindent Theorem~\ref{thm:hms} in the previous section can also be generalized to this situation for $\cO^\omega$ and $\cO_{M^\chk(U)}^\hol$. This result is summarized in the following theorem. Its proof is again to use spectral sequences associated to energy filtrations to reduce to classical results. In this case instead of using classical Koszul duality we use classical Fourier transform for families, see for instance~\cite{BMP}. We shall not repeat the proof here.

\begin{theorem}
\label{thm:hms2}
The composition of $A_\infty$ functors
\[\Fuk^\pi(M)\stackrel{P}{\ra} \tw(\cO^\omega) \stackrel{\Phi^\cP}{\ra} \Tw(\cO_{M^\chk(U)}^\hol)\]
is a quasi-equivalence onto its image. Here the first functor $P$ was defined in Theorem~\ref{thm:lag}.
\end{theorem}

\remark Here we need to use $\Tw(\cO_{M^\chk(U)}^\hol)$ to include infinite rank objects over $\cO_{M^\chk(U)}^\hol$. It is an interesting question to do homological perturbation on $\Phi^\cP\circ P(\sL_p(\alpha))$ to reduce to an object of finite rank. This problem might be related to the appearance of theta functions in mirror symmetry~\footnote{The assertion is motivated from the case of elliptic curves.}.

\section{Homological mirror symmetry on toric manifolds}
\label{sec:functor}

As an immediate application of our general theory we prove a version of homological mirror symmetry between a toric symplectic manifold and its Landau-Ginzburg mirror. 

\begin{theorem}
\label{thm:toric}
Let $M$ be a compact smooth toric symplectic manifold, and denote by $\pi:M(\Delta^\interior)\ra \Delta^\interior$ the Lagrangian torus fibration over the interior of the polytope of $M$. Then there exists an $A_\infty$ functor $\Psi: \Fuk^\pi(M) \ra \tw(\cO^\hol_{T\Delta^\interior})$ which is a quasi-equivalence onto its image.
\end{theorem}

\proof We take $J$ to be the canonical integrable complex structure on $M$. Since biholomorphisms of $M$ acts transitively on torus fibers over $\Delta^\interior$, we get a sheaf of $A_\infty$ algebras over $\Delta^\interior$ by constructions in Section~\ref{sec:symp}. Moreover it is known after~\cite{FOOO} that the weak unobstructedness assumption~\ref{ass:wua} holds in this context. Thus all results in this paper applies to this situation, and the theorem is simply an example of Theorem~\ref{thm:hms}.

\medskip
\remark The functor $\Psi$ in the above theorem is simply the composition $\Phi^\tau\circ P$ where $P$ is the propagation functor used in Theorem~\ref{thm:lag}, and $\Phi^\tau$ is the Koszul duality functor $\Phi^\tau$ defined in Section~\ref{sec:koszul}. The image of $\Psi$ consists of matrix factorizations of $W$ on $T\Delta^\interior$ which can be calculated by the formula
\[ Q (f\otimes e_I):=\sum_{k\geq 0,l\geq 0} m_{k+l+1}(\tau^l,f\otimes e_I,\theta^k).\]
It is plausible that these objects split generate $\tw(\cO^\hol_{T\Delta^\interior})$. But we do not know how to prove this generation result in general. Some of the difficulties are
\begin{itemize}
\item working over Novikov ring instead of $\C$;
\item the map $Q$ mixes various wedge degrees.
\end{itemize}
In the Fano case we can specialize to $T=e^{-1}$ and work over $\C$, which enables us to get around the first issue. When the dimension is less than or equal to two the inhomogeneity does not appear, which allows us to prove the following.

\begin{theorem}
\label{thm:toricfano}
Let $M$ be a compact smooth toric Fano symplectic manifold of dimension less or equal to two. In this case we can work over $\C$ by evaluating the parameter $T$ at $e^{-1}$. Then there is a functor $\Psi^\C:\Fuk^\pi(M,\C) \ra \tw(\cO^\hol_{T\Delta^\interior}\otimes\C)$ which is a quasi-equivalence of $A_\infty$ categories.
\end{theorem}

\proof The functor $\Psi^\C$ is simply the reduction of $\Psi$ at the evaluation $T=e^{-1}$, which is valid under the Fano condition. By the previous theorem it suffices to show that the image of $\Psi^\C$ split generates the target category $\tw(\cO^\hol_{T\Delta^\interior}\otimes\C)$. This generation follows from explicitly computing the operator $Q$ and using the generation result of T. Dyckerhoff~\cite{Dyck} Section $4$. Note that this generation result requires $W$ to have isolated singularities which was proved in~\cite{FOOOtoric} Theorem $10.4$.

Next we compute the operator $Q$. In the following the degree $|I|$ of $e_I$ is referred to as wedge degree. Recall the operator $Q_\beta$ is of wedge degree $1-\mu(\beta)$. Since for toric manifolds there are no negative Maslov index holomorphic disks, the operator $Q=\sum_{\beta\in G} Q_\beta$ has only one part $Q_0$ that increases the wedge degree. As we saw in the proof of Proposition~\ref{prop:torus} $Q_0$ is the Koszul differential associated to the regular sequence $z_i -p_i-\sqrt{-1}\alpha_i$.

If the dimension is less than or equal to two, then the operator $Q_\beta$ is necessary of wedge degree $-1$ corresponding to $\mu(\beta)=2$. Such type of matrix factorizations is shown to split generate $\tw(\cO^\hol_{T\Delta^\interior}\otimes\C)$ by~\cite{Dyck} Section $4$. The theorem is proved.

\paragraph{An example: $\C P^1$.} A particularly simple example is the case $M=\C P^1$. Since it is Fano we shall work over $\C$. With appropriate choice of its symplectic form we assume $U=(0,1)\subset \R$ as is in~\cite{FOOOtoric} Section $5$. Let $e$ be a trivialization of $R^1\pi_*\Z$, and let $x$, $y^\chk$, $y$ be associated affine coordinates. It is known that the $A_\infty$ algebra associated to the Lagrangian torus fiber $L_{x}$ for $x\in (0,1)$ is a two dimensional vector space generated by $\bone, e$ with $\bone$ a strict unit. All the rest $A_\infty$ products are
\begin{align*}
m_0 &= \exp(-x)+ \exp(x-1);\\
m_1 (e) &= \exp(-x)-\exp(x-1); \\
&\cdots;\\
m_k (e^{\otimes k}) &=\frac{1}{k!}[ \exp(-x)+ (-1)^k \exp(x-1)];\\
&\cdots.
\end{align*}
The potential function is equal to
\begin{align*}
W(x,-\sqrt{-1}y)&= \sum_{i=0}^\infty m_k( (-\sqrt{-1}y e)^{\otimes k})\\
&= \sum_{i=0}^\infty \frac{1}{k!} (-\sqrt{-1}y)^k [ \exp(-x)+ (-1)^k \exp(x-1)]\\
&=\exp(-z)+\exp(z-1).
\end{align*}
Thus $W$ is indeed a holomorphic function. Let $a\in \R$ be a real number, and let $u\in (0,1)$ be a point. Consider the Lagrangian brane $(L_u,-\sqrt{-1}a)$. From this data we get define an $A_\infty$ module $\sL_u(\sqrt{-1}a)$ over $\cO_{M(U)}^{\omega,\can}$ with internal curvature $W(u,a)$ by constructions in Section~\ref{sec:lag}. Let us describe its image under the Koszul functor $\Phi^\tau$. It suffice to compute the operator $Q$ on generators $\bone$ and $e$. For this we have
\begin{align*}
Q(1)&=\sum_{k,l}  m_{k+l+1}(\tau^l,\bone,\theta^k)\\
&=\sum_{k,l} (x-u-\sqrt{-1}a)^k(-\sqrt{-1}y)^l m_{k+l+1}(e^{\otimes l},\bone,e^{\otimes k})\\
&= e\otimes [(x-u-\sqrt{-1}a)-(-\sqrt{-1}y)] \mbox{\;\;\; (by our sign convention)}\\
 &=e\otimes [z-u-\sqrt{-1}a].
\end{align*}
A more technical computation of $Q(e)$ by formula gives $\frac{W-W(u,a)}{u+\sqrt{-1}a-z}$, and hence $Q^2+[W-W(u,a)]\id=0$ as is expected from the general theory. Thus we see that $\Phi^\tau(\sL_u(\sqrt{-1}a))$ is a matrix factorization of $-[W-W(u,a)]$.

\appendix

\section{Modules with internal curvature.}
\label{app:modules}

Let $A$ be an $A_\infty$ algebra. In this section we define $A_\infty$ modules over $A$ possibly with an internal curvature. We show how weak Maurer-Cartan elements of $A$ give rise to such structures. These algebraic constructions naturally occur in Lagrangian Floer theory.

Throughout the construction we work over a base ring $R$, and all modules considered here are free $R$-modules. We follow the sign convention used in~\cite{FOOO}. We refer to~\cite{Keller} Section $3$ and $4$ for basics of $A_\infty$ algebras, homomorphisms and modules.

\paragraph{$A_\infty$ modules.} An $A_\infty$ module $M$ over an $A_\infty$ algebra $A$ is defined by a collection of maps $\rho_k(-;-): (A^{\otimes k})\otimes M\ra M$ of degree $1-k$ such that
\begin{equation}
\label{eq:module}
 \sum_{i+j=N} \rho_i(\id^i;\rho_j(\id^j;-))+ \sum_{r+s+t=N} \rho_{r+t+1} (\id^r,m_s(\id^s),\id^t;-)=0
\end{equation}
for all $N\geq 0$. When applied to elements $(a_1\otimes\cdots\otimes a_N\otimes x)\in A^{\otimes N}\otimes M$ extra signs come out by Koszul sign rule, for example when $N=0,1$ the above relation reads

\begin{align*}
\rho_0(\rho_0(x))&+\rho_1(m_0;x)=0;\\
\rho_0(\rho_1(a;x))&+(-1)^{|a|-1}\rho_1(a;\rho_0(x))+\\
+&\rho_{1}(m_1(a);x)+\rho_{2}(m_0,a;x)+(-1)^{|a|-1}\rho_2(a,m_0;x)=0.
\end{align*}

\noindent Using the Bar construction we can interpret an $A_\infty$ module structure on a $R$-module $M$ as an $A_\infty$ homomorphism $\rho:A\ra \End(M)$~\footnote{Note that the product on the graded matrix algebra $\End(M)$ is defined by $(\phi\otimes\psi)\mapsto(-1)^{|\phi|}\phi\circ\psi$
where the sign appears due to our sign convention.}. Recall an $A_\infty$ homomorphism between This correspondence is explicitly given by 
\[\left\{\rho_k\right\}_{k=0}^{\infty} \mapsto \rho:=\prod_{k=0}^\infty \rho_k\in \Hom_{A_\infty}(A,\End(M))\]

\paragraph{$A_\infty$-modules with internal curvature.} For applications in this paper we need to introduce a weaker notion of modules: those endowed with ``internal curvatures". For its definition we will fix $\lambda\in R$ an even element in the ground ring.
\begin{definition}
\label{def:module}
An $A_\infty$ module $M$ over $A$ with internal curvature $\lambda$ is defined by structure maps $\rho_k(-;-): (A^{\otimes k})\otimes M\ra M$ of degree $1-k$ which satisfies the same axioms as in equation~\ref{eq:module} except for $N=0$ in which case we require that
\[\rho_0(\rho_0(x))+\rho_1(m_0;x)=\lambda\id_M.\]
\end{definition}

\noindent From the point of view of $A_\infty$ homomorphisms, we can add the element $\lambda\id_M$ as a curvature element for the matrix algebra $\End(M)$, and we denote the resulting curved algebra by $\End^\lambda(M)$. Then straight-forward computation shows that an $A_\infty$ module $M$ over $A$ with internal curvature $\lambda$ is the same as an $A_\infty$ homomorphism $\rho:A\ra \End^\lambda(M)$.

\paragraph{From weak Maurer-Cartan elements to modules with internal curvature.} Let us see how $A_\infty$ modules with internal curvature can arise from a weak Maurer-Cartan element of $A$. Recall if $A$ is an $A_\infty$ algebra with a strict unit $\bone$, an odd element $b\in A^1$ is a weak Maurer-Cartan element if we have
\[ \sum_{k=0}^\infty m_k(b,\cdots,b)=\lambda\bone\]
for some even element $\lambda\in R$. Using such an element we can define an $A_\infty$ module structure on the same underlying space of $A$ which has internal curvature $\lambda$. We denote this $A_\infty$ module by $A^b$. Its structure maps are defined by

\[\rho_k^b(a_1,\cdots,a_k;x) :=\sum_{i=0}^{\infty} m_{i+k+1}(a_1,\cdots,a_k,x,b^{\otimes i})\]

\noindent Let us check the first axiom, i.e. $\rho^b_0(\rho^b_0(x))+\rho^b_1(m_0;x)=\lambda\id_{A^b}$.

\begin{align*}
&\rho^b_0(\rho^b_0(x)) =\sum_{i,j} m_{i+j+1}(m_{i+1}(x,b^{\otimes i}),b^{\otimes j}) \\
&=-\sum_{r\geq 0, s\geq 0, t\geq 0} (-1)^{|x|-1} m_{r+t+2} (x,b^{\otimes r},m_s(b^{\otimes s}),b^{\otimes t}) -\sum_{k\geq 0} m_{k+2}(m_0,x,b^{\otimes k}) \\
&=(-1)^{|x|}m(x,\lambda\bone)-\rho_1^b(m_0;x) \mbox{\;\;\; (by the weak Mauer-Cartan equation)}\\
&=\lambda\id -\rho_1^b(m_0;x) \mbox{\;\;\; (by our sign convection of strict unit).}
\end{align*}

\noindent The rest identities can be checked similarly using the fact that $\lambda\bone$ is a multiple of strict unit, and hence does not contribute to higher products. We denote by $\rho^b:A\ra \End^\lambda (A^b)$ the corresponding $A_\infty$ homomorphism.

\paragraph{Twisted complexes.} The category of $A_\infty$ modules are usually defined as a differential graded category. However for purposes of the current paper it is better to use the category of twisted complexes of $A$ which is an $A_\infty$ category. We refer to the paper of B. Keller~\cite{Keller} Section $7$ for details of these categorical constructions. The category of twisted complexes over an $A_\infty$ algebra will be denoted by $\tw(A)$. Intuitively this can be thought of as the $A_\infty$ analogue of differential graded modules over an algebra $A$ that are free of finite rank.

To include modules with internal curvatures we need to modify slightly the definition of $\tw(A)$. In the following we explain this modification. This modified version of $\tw(A)$ is a direct sum of $R$-linear categories
\[ \tw(A):=\coprod_{\lambda\in R^\even} \tw^\lambda(A).\]
For each $\lambda\in R^\even$ the category $\tw^\lambda(A)$ consists of twisted complexes over $A$ with internal curvature $\lambda$. Thus the conventional definition of $\tw(A)$ corresponds to $\tw^0(A)$ in our notation. Let us explain in more detail the construction of $\tw^\lambda(A)$.

For each $\lambda$ the objects of $\tw^\lambda(A)$ are pairs $(V,b)$ where $V$ is a finite rank $\Z/2\Z$-graded free $R$-module, and $b$ is a weak Maurer-Cartan element of the tensor product $A\otimes \End_R(V)$ with internal curvature $\lambda$. By constructions in the previous paragraph these data give rise to an $A_\infty$ module over $A\otimes V$ with internal curvature $\lambda$. Strictly speaking in the previous paragraph we only dealt with the case when $V$ is of rank one over $R$, but the general case only requires more index.

Let us illustrate the morphism space between two pairs $(V,b)$ and $(W,\delta)$ when both $V$ and $W$ is one dimensional. In this case the $\Hom$ space, as a graded $R$-module, is simply $A$ itself. It is endowed with a differential $d$ twisted by $b$ and $\delta$ which is explicitly given by formula
\[ a\mapsto \sum_{k,l} m_{k+l+1}(b^k,a,\delta^l).\]
Using $A_\infty$ relations and Maurer-Cartan equations one shows that $d^2=[F(b)-F(\delta)]\id=[\lambda-\lambda]\id=0$ where $F(b)$ and $F(\delta)$ are internal curvatures of $b$ and $\delta$. This explains the reason why twisted complexes with different internal curvatures do not interact with each other. For the general case when $V$ and $W$ are of any finite rank, the definition is similar using matrix compositions. We refer the details to~\cite{Keller} Section $8$.

\paragraph{Upper-triangular condition.} Finally we end this appendix with an important technical point involved in the construction of twisted complexes. In conventional definitions one usually assumes that the (weak) Maurer-Cartan element $b\in A\otimes \End_R(V)$ to satisfy strict upper-triangular condition. This has two important implications. Namely this assumption implies convergence of Maurer-Cartan equation and also the convergence of twisted differential. Secondly it also implies that the homotopy category of $\tw(A)$ embeds fully into the derived category of $A$; moreover the image of this embedding is simply the triangulated closure of $A$ as an $A_\infty$ module over itself (this works when $m_0$ of $A$ is trivial).

While working with such a condition has nice homological implications, it is too restrictive for applications in mirror symmetry. Indeed it follows from the upper-triangularity that the only rank one twisted complex is $A$ itself if $m_0$ vanishes. But as is shown in Section~\ref{sec:lag} we would like to associate to each Lagrangian torus fiber a non-trivial rank one twisted complex. Thus we would like to work with twisted complexes which might not satisfy the upper-triangular condition. In this case convergence of relevant series is not automatic, and needs to be taken care of separately. For Lagrangian Floer theory as needed in this paper this convergence follows from results in~\cite{Fukaya}. Secondly the homotopy category of $\tw(A)$ in our definition might not admit a fully faithfully embedding into the derived category of $A$.

\section{Koszul duality as Fourier-Mukai transform}
\label{app:fm}
In this appendix we construct a Koszul duality functor as a type of affine version of Fourier-Mukai transform. We also define such a functor on modules with internal curvatures. It follows from our definition the Koszul functors preserve internal curvatures of modules.

We will need to use another sign convection~\cite{Keller} since the previous sign convention is not convenient to deal with tensor product of algebras. To avoid possible confusions from using two different signs, we first clarify the relationship between them. 

\paragraph{Sign conventions.} In the sign convention used in the previous appendix, the maps $m_k$ are considered as degree one maps between suspensions $(A[1])^{\otimes k}\ra A[1]$. The advantage of doing so is that there is no signs in the $A_\infty$ algebra axioms, i.e. for each $n\geq 0$ we have
\[ \sum_{r+s+t=n} m_{r+1+t}(\id^r\otimes m_s \otimes \id^t)=0.\]
When applying to elements we get signs by the Koszul rule. Using maps $m_k$ we can define its corresponding linear maps $m_k^\epsilon:A^{\otimes k}\ra A$ by requiring the following diagram to be commutative:
\[\begin{CD}
A^{\otimes k} @>m_k^\epsilon>> A \\
@VV [1]^{\otimes k} V                   @VV [1] V \\
(A[1])^{\otimes k} @>m_k>>  A[1].
\end{CD}\]
Here the map $[1]: A\ra A[1]$ defined by $a\mapsto a[1]$ is the identity map on the underlying $R$-module, but since it is a degree one map it yields signs when applied to tensor products by Koszul sign rule. Explicitly we have
\[ m_k^\epsilon(a_1,\cdots,a_k)[1]=(-1)^{\epsilon_k} m_k(a_1[1],\cdots,a_k[1])\]
where $\epsilon_k=\sum_{i=1}^k |a_i|(k-i)$. The above identity applied to the $A_\infty$ axioms of $m_k$ yields
\[(-1)^{\epsilon_n}\sum_{r+s+t=n} (-1)^{r+st} m^\epsilon_{r+1+t}(\id^r\otimes m^\epsilon_s \otimes \id^t)=0.\]
Dividing the sign $(-1)^{\epsilon_n}$ gives the $A_\infty$ axioms for the maps $m^\epsilon_k$. Using this relationship between $m_k$ and $m_k^\epsilon$ we can freely pass from one to the other. For instance the weak $\MC$ equation expressed using $m_k^\epsilon$ reads
\[\sum_{k=0}^\infty (-1)^{\frac{k(k-1)}{2}}m_k^\epsilon(b^{\otimes k}) =\lambda\bone.\]
A strict unit $\bone$ in the $\epsilon$-sign convention becomes
\begin{align*}
m_2^\epsilon(1,x)=m_2^\epsilon(x,1)&=x\mbox{\;\;\; and} \\
m_k^\epsilon(a_1,\cdots,a_i,\bone,\cdots,a_{k-1})&=0 \mbox{\;\;\; for all $k\neq 2$.}
\end{align*}
\paragraph{Tensor product.} Let $A$ be a strict unital $A_\infty$ algebra, and let $B$ be a curved differential graded algebra. Form their tensor product $B\otimes A$ which is an $A_\infty$ algebra with structure maps defined by
\begin{align*}
m_0^\epsilon&:=\bone\otimes m_0^\epsilon +W\otimes\bone;\\
m_1^\epsilon(b\otimes a) &:=db\otimes a+(-1)^{|b|} b\otimes m_1^\epsilon (a);\\
m_k^\epsilon(b_1\otimes a_1,\cdots,b_k\otimes a_k)&:=(-1)^{\eta_k}(b_1\cdots b_k)\otimes m_k^\epsilon(a_1,\cdots,a_k)  \mbox{\;\;\; for $k\geq 2$.}
\end{align*}
Here $W$ is the curvature term of $B$, the two $\bone$'s are units, and the sign in the last equation is $\eta_k=\sum_{i=1}^{k-1} |a_i|(|b_{i+1}|+\cdots+|b_k|)$. We have abused the notation $m_k^\epsilon$ for both structure maps on $A$ and $B\otimes A$. Since they are applied to different types of elements, no confusion can arise by doing so.

\paragraph{Koszul duality as Fourier-Mukai transform.} Next we describe a construction of a functor $\Phi^\tau:\tw(A)\ra\tw(B)$ associated to a given Maurer-Cartan element $\tau\in B\otimes A$. For simplicity we will assume that $A$ is of finite rank over $R$. This is for the purpose that $\Phi^\tau$ lands inside $\tw(B)$, i.e. it is of finite rank over $B$. If we replace the target category by $\Tw(B)$ consisting of twisted complexes over $B$ of possibly infinite rank then all constructions below still go through.

The intuitive idea to construct such a functor is that the element $\tau$ determines an $A_\infty$ module $(B\otimes A)^\tau$ over $B\otimes A$ which can be viewed as a kernel for an integral transform from $\tw(A)$ to $\tw(B)$. To realize this idea we proceed as follows.

Given an $A_\infty$ module $M$ with internal curvature $\lambda$, we denote by $\rho^M$ the corresponding $A_\infty$ homomorphism $A\ra \End^\lambda(M)$. For $M\in\tw(A)$ by our assumption that $A$ is of finite rank, it follows that $M$ is also of finite rank. For general $M$ we assume $M$ is of finite rank below.

The map $\rho^M$ induces  another $A_\infty$ homomorphism $\rho^M_B: B\otimes A\ra B\otimes\End^\lambda(M)$ by scalar extension to $B$. Using $\rho^M_B$ we can push forward the given Maurer-Cartan element $\tau$ to get a Maurer-Cartan element of $\End^\lambda(M)\otimes B$. Such a Maurer-Cartan element by definition is a twisted complex structure on $B\otimes_R M$ with internal curvature $\lambda$. The following Theorem gives a more explicit description of this construction with formulas.

\begin{theorem}
\label{thm:duality}
The maps $(\rho^M_B)_k:(B\otimes A)^k \ra B\otimes\End^\lambda(M)$ defined by
\begin{align*}
(\rho^M_B)_0&:=\bone\otimes\rho^M_0;\\
(\rho^M_B)_1(b_1\otimes a_1)&:=b_1\otimes\rho^M_1(a_1);\\
(\rho^M_B)_k(b_1\otimes a_1,\cdots,b_k\otimes a_k) &:=(-1)^{\eta_k} (b_1\cdots b_k)\otimes \rho^M_k(a_1,\cdots,a_k)
\end{align*}
form an $A_\infty$ homomorphism $\rho^M_B: B\otimes A \ra B\otimes\End^\lambda(M)$. Moreover if $\tau\in B\otimes A$ is a Maurer-Cartan element, its push-forward $Q:=(\rho_B^M)_*\tau=\sum_{k=0}^\infty (-1)^{\frac{k(k-1)}{2}} (\rho^M_B)_k(\tau^{\otimes k})$ is a Maurer-Cartan element of $B\otimes \End^\lambda(M)$, i.e. we have 
\[(\lambda-W)\id+[d,Q]-Q^2=0.\]
Here $-W$ is the curvature of $B$, and $d$ is its differential.
\end{theorem}

\proof The proof is straightforward verifications of formulas and keeping track of signs. We omit it here.

\medskip
\noindent By the above theorem we described what $\Phi^\tau$ does on the level of objects. Namely we define $\Phi^\tau(M)$ to be the twisted complex on $B\otimes_R M$ defined by the weak Maurer-Cartan element $(\rho_B^M)_*\tau$. Note that $\Phi^\tau(M)$ and $M$ have the same internal curvature, i.e. we have a map
\[ \Phi^\tau: \tw^\lambda(A) \ra \tw^\lambda(B)\]
for each $\lambda\in R^\even$. Next we describe $\Phi^\tau$ on the level of morphisms. We can define $\Phi^\tau$ as an $A_\infty$ functor from $\tw(A)$ to $\tw(B)$. Let us illustrate this for a rank one twisted complex $A^b$ over $A$ with internal curvature $\lambda$. The space $\End(A^b)$ is an $A_\infty$ algebra with structure maps $m_k^b$ defined by 
\[m_k^b(a_1,\cdots,a_k):=\sum_{i_0\geq 0,\cdots,i_k\geq 0} m_{i_0+i_1+\cdots+i_k+k} (b^{i_0},a_1,b^{i_1},\cdots,b^{i_{k-1}},a_k,b^{i_k}).\]
We need to define an $A_\infty$ homomorphism from $\End(A^b)$ to the differential graded algebra $\End(\Phi^\tau(A^b))$. Since $\Phi^\tau(A^b)=B\otimes_R A$ as a $R$-module, we use structure maps $m_k$ on the tensor product $B\otimes A$ to describe this homomorphism. Explicitly it is given by
\[ \Phi^\tau(a_1,\cdots,a_k)(x):= \sum_{l\geq 0,i_0\geq 0,\cdots,i_k\geq 0} m_{l+k+1+i_0+\cdots+i_k}(\tau^l,x,b^{i_0},a_1,b^{i_1},\cdots,a_k,b^{i_k}).\]
The case of higher rank twisted complexes requires no more than putting more index into the above equation. We refer to Section $7.3$ of K. Lef\`evre-Hasegawa's thesis~\cite{LH} for a more detailed discussion of this $A_\infty$ homomorphism.

\end{document}